\documentclass[12pt]{article}
\usepackage{latexsym,epsfig,stmaryrd}
\usepackage{amssymb}

\title{Measure Functions for Frames}

\author{\bf Radu Balan \\  \small Siemens Corporate Research \\ 
\small Princeton, NJ 08540 \\ \bf Zeph Landau \\  \small The City College of New York \\ \small New York, NY 10031}

\newcommand{\qed}{\mbox{$\Box$}}

\newcommand{\cs}{X}
\newcommand{\rs}{X^{\R}}
\newcommand{\cse}{\rs}

\newcommand{\cxs}{X^{\C}}

\newcommand{\plim}{{p\_\lim } }

\newcommand{\ec}{{\mathcal E}}

\newcommand{\G}{{\mathcal G}}
\newcommand{\R}{{\mathbf R}}
\newcommand{\Z}{{\mathbf Z}}
\newcommand{\C}{{\mathbf C}}
\newcommand{\N}{{\mathbf N}}

\newcommand{\pf}{{\it Proof: }}
\newcommand{\eps}{{\varepsilon}}
\newcommand{\ip}[2]{\langle#1,#2\rangle}
\newcommand{\norm}[1]{\|#1\|}
\newcommand{\ltwoi}{l^2(I)}

\newcommand{\tFc}{{\tilde{\mathcal{F}}}}
\newcommand{\tf}{{\tilde{f}}}
\newcommand{\fc}{{\cal F}}
\newcommand{\fci}{{\cal F}[I]}
\newcommand{\gc}{{\cal G}}
\newcommand{\cc}{{\cal C}}
\newcommand{\hc}{{\cal H}}

\newcommand{\tm}{{\tilde{m}}}

\newcommand{\simm}{\approx}
\renewcommand{\leqq}{{\trianglelefteqslant}}

\newcommand{\accc}{p}
\newcommand{\E}{{\bf E}}

\newcommand{\xv}{{\bf x}}
\newcommand{\yv}{{\bf y}}
\newcommand{\zv}{{\bf z}}
\newcommand{\bv}{{\bf b}}

\newcommand{\iv}{{\bf i}}

\newcommand{\xvt}{{\tilde{\bf x}}}
\newcommand{\cstarm}{\cc^*(M)}
\newcommand{\cstarw}{\cc^*(W)}
\newcommand{\cstarwcx}{\cc^*_{\C}(W)}
\newcommand{\cstarn}{\cc^*(\N^*)}

\newcommand{\cstarno}{\cc^*(\N^0)}
\newcommand{\cstarwo}{\cc^*(W^0)}

\newcommand{\Cstar}{\mbox{$C^*$}}

\newcommand{\Bltwoi}{\mbox{$B(\ltwoi)$}}

\newcommand{\mt}{{\tilde{m}}}
\newcommand{\metricindexset}{uniform metric index set}

\newcommand{\No}{{\N^0}}
\newcommand{\set}[1]{\{#1\}}
\newcommand{\fperp}{perpendicular-normal }
\newcommand{\mbar}{\overline{m}}
\newcommand{\ignore}[1]{}
\newcommand{\al}{\alpha}
\newcommand{\be}{\beta}

\renewcommand{\span}[1]{\overline{\mbox{span}\{#1\}}}

\newtheorem{Theorem}{Theorem}[section]
\newtheorem{Corollary}[Theorem]{Corollary}
\newtheorem{Lemma}[Theorem]{Lemma}
\newtheorem{Proposition}[Theorem]{Proposition}
\newtheorem{Definition}[Theorem]{Definition}
\newtheorem{Remark}[Theorem]{Remark}
\newtheorem{Example}[Theorem]{Example}

\begin{document}

\maketitle


\begin{abstract}
This paper addresses the natural question: ``How should frames be
compared?''  We answer this question by quantifying the
overcompleteness of all frames with the same index set.  We introduce
the concept of a {\it frame measure function}: a function which maps
each frame to a continuous function.  The comparison of these
functions induces an equivalence and partial order that allows for a
meaningful comparison of frames indexed by the same set.  We define
the {\it ultrafilter measure function}, an explicit frame measure
function that we show is contained both algebraically and
topologically inside all frame measure functions.  We explore
additional properties of frame measure functions, showing that they
are additive on a large class of supersets-- those that come from so
called {\it non-expansive} frames.  We apply our results to the Gabor
setting, computing the frame measure function of Gabor frames and
establishing a new result about supersets of Gabor
frames. 
\end{abstract}


\section{Introduction}\label{sec1}

Let $H$ be a separable Hilbert space and $I$ a countable index set.
A sequence $\fc = \set{f_i}_{i \in I}$ of elements of $H$ is a
{\em frame} for $H$ if there exist constants $A$, $B>0$ such that
\begin{equation}
\label{framedef}
\forall\, h \in H, \quad
A \, \norm{h}^2 \le \sum_{i \in I} |\ip{h}{f_i}|^2 \le B \, \norm{h}^2.
\end{equation}
The numbers $A$, $B$ are called {\em lower} and {\em upper frame bounds},
respectively.
 Frames were first introduced by Duffin and Schaeffer \cite{dusc52}
in the context of nonharmonic Fourier series, and today frames
play important roles in many applications in mathematics, science,
and engineering.
We refer to the monograph of Daubechies \cite{da92} or the
research-tutorial \cite{Cas00} for basic properties of frames.

Central, both theoretically and practically, to the interest in frames
has been their overcomplete nature; the strength of this
overcompleteness is the ability of a frame to express arbitrary
vectors as a linear combination in a ``redundant'' way.  Until
recently, for infinite dimensional frames, the
overcompleteness or redundancy has only referred to a qualitative
feature of frames.  A notable exception is in the case of Gabor frames
where many works have connected essential features of the frames to
quantities associated to the density of the associated lattice of time
and frequency shifts (\cite{he06-1} and references therein).  
Recently, the work in
\cite{bacahela03,bacahela03-1,bacahela06,bacahela06-1} examined and explored
the notion of {\it excess} of a frame, i.e. the maximal number of
frame elements that could be removed while keeping the remaining
elements a frame for the same span.  A quantitative approach to
certain frames with infinite excess was given in \cite{bacahela06,bacahela06-1}
which introduced a general notion of a localized frame and, among
other results, provided nice quantitative measures associated to this
class of frames.

This paper addresses the natural question: ``How should frames be
compared?''  We answer this question by quantifying the
overcompleteness of all frames with the same index set.  We describe a
new equivalence relation and partial order on these frames.  We
introduce the central tool for working with this partial order: the
{\it frame measure function } which maps each frame to a continous
function.  The frame measure functions are compatible with our
equivalence relation, namely two frames are equivalent if and
only if their frame measure functions are equal (pointwise as
continuous functions) and one frame dominates another if their
frame measure functions have the corresponding dominance
(pointwise).  This results in a quantification of frames that reflects
the partial order and leads to a meaningful quantitative definition of
the overcompleteness of a frame.

Though equivalence of frames with an infinite number of elements has
been considered previously (see \cite{ba99-2,dailarson98}) and a standard
notion of equivalence for frames exist, the size of each equivalence
class is too small; it is fundamentally unsatisfying as it
distinguishes frames, that from a signal processing point of view, are
equivalent.

In contrast, the equivalence relation, partial order and frame measure
function introduced here have the following desirable properties (that
are not present in the standard equivalence relation):

\begin{itemize}
\item  The equivalence relation groups together all Riesz bases.
\item The equivalence relation groups together all frames that differ
by a finite permutation of their elements or by arbitrary phase change
of their elements.
\item From an information theory point of view, the equivalence
relation groups together frames that transmit signals with similar
variances due to noise.
\item The values of the measure function are linked to the amount of
excess of the frame.
\item For a large class of frames (those that are called {\it
non-expansive}) any frame measure function is additive on supersets,
namely, the frame measure function applied to the frame $\{f_i \oplus
g_i\}_{i\in I}$ acting on $H_1 \oplus H_2$ is equal to the sum of the
frame measure function applied to the two frames $\{f_i\}_{i \in I}$
acting on $H_1$ and $\{g_i\}_{i \in I}$ acting on $H_2$.
\item The values of the frame measure function for Gabor frames
are shown to correspond to the density in the time-frequency plane of
the shifts associated to the frame.
\end{itemize}

The focus of this work is to explore the properties of the equivalence
relation, partial order, and frame measure functions.  In addition to
showing the above listed facts, we describe a specific frame measure
function, the {\it ultrafilter frame measure function }-- a function
from the set of all frames indexed by a set $I$ (denoted by $\fci$) to
the set of continuous functions on the compact space consisting of the
free ultrafilters.  We show that every frame measure function contains
a copy of the ultrafilter frame measure function.  In addition, as
with representation theory, we define separable, reducible, and
minimal frame measure functions and show that all minimal frame
measure functions are topologically equivalent to the ultrafilter
frame measure function.

We apply this theory to the Gabor setting.  In addition to
computing the measure of Gabor frames, we apply our results to Gabor
supersets, showing new necessary conditions on the densities of the
time-frequency shifts of the individual Gabor frames.

Finally we propose that the reciprocal of the measure function be
defined to be the redundancy for an infinite frame.  Redundancy, an
often referred to qualitative feature of frames, has eluded a
meaningful quantitative definition for infinite frames.  Using the
results of this work, we justify our definition of redundancy by both
showing it to be quantitatively meaningful and a natural
generalization of redundancy for finite frames.

A striking feature of these ideas is the variety of mathematical areas
that are involved.  The fundamental objects, frames, are objects of
considerable interest to the signal processing community.  The
motivation for our definitions of frame equivalence and comparison
come from both information theoretic and operator theoretic
considerations.  The ideas and tools that drive the results are mainly
operator theoretic and topological.

The equivalence relation, partial order, and frame measure
functions introduced here are a function of certain averages of the
terms $\ip{f_i}{\tilde{f}_i}$ of a given frame $\{f_i\}_{i\in I}$ (where
$\{\tilde{f}_i\}_{i\in I}$ is the canonical dual frame to
$\{f_I\}_{I\in I}$).  These are the same averages that play a central
role in the two papers \cite{bacahela06,bacahela06-1} which introduce the notion of
localized frames.  In this work, our goal is to compare all frames
that are indexed by the same fixed index set but which possibly lie in
different Hilbert spaces; we require no special localized structure
for the frames.  In contrast, in \cite{bacahela06,bacahela06-1} the situation
considered is that of frames which all lie in the same Hilbert space
that are indexed by different sets.  An index set map is introduced
and when this index map is chosen so that the frame is localized,
powerful results are obtained relating a feature of the index map
(density), to certain averages of $\ip{f_i}{\tilde{f}_i}$ (relative
measure).  Despite the differences in approach between
\cite{bacahela06,bacahela06-1} and this work, there is significant intersection and
interelation of ideas.  Specifically, where the settings are
compatible, the notion of a {\it non-expansive frame} introduced here
is the same as the notion of a $l^2$ localized frame of
\cite{bacahela06,bacahela06-1}.  In addition, we use specific results of
\cite{bacahela06-1} to compute the ultrafilter frame measure function of
Gabor frames.

The work is organized as follows.  The equivalence relation and
partial order is introduced and initially explored in Section
\ref{s:equiv}.  Section \ref{s:ult} defines and proves essential
properties of the {\it ultrafilter frame measure function}.  The
general notion of a frame measure function is defined and core
properties are proven in Section \ref{s:gen}.  Of particular note is
Corollary \ref{c:frameembed} which shows that every frame measure
function contains an algebraic copy of the ultrafilter frame measure
function.  Section \ref{s:top} examines the topological properties of
the frame measure function, showing, among other things, that in a
certain sense, that the ultrafilter frame measure function is the
unique minimal frame measure function. We extend the frame measure
functions ideas to the space of operators in Section \ref{s:9} and
introduce the core concept of a non-expansive operator.  Section
\ref{sec4.2} applies these ideas to supersets to prove Theorem
\ref{t:redundancy} which establishes that frame measure functions are
additive on superframes comprised of non-expansive frames.  Section
\ref{s:index} examines the connection between the measure function and
the index set.  Section \ref{s:gab} applies the results to the Gabor
setting, computing the frame measure function of Gabor frames and
establishing a new result about supersets of Gabor frames.  Finally,
section \ref{s:red} defines and explores the properties of the
redundancy function for infinite frames.  The Appendices cover some
background material on supersets and ultrafilters.

\section{Notation and Preliminaries}\label{s:2}

\subsection{Basic Notation} \label{s:2.1}
For any set $S$, $|S|$ will denote the number of elements in $S$.
Throughout this paper $I$ will be a fixed countable index set
accompanied by a decomposition into a nested union (indexed by the
positive integers $1,2,\dots$) of finite subsets.  That is,
\begin{eqnarray}
 I_1\subset I_2 \subset \cdots \subset & I_n & \subset I_{n+1} \subset \cdots
\subset I  \label{eq.e.1.2.1aa} \\
 |I_n| & < & \infty  \label{eq.e.1.2.2aa} \\
 \cup_{n\geq 1} I_n & = & I \label{eq.e.1.2.3aa}.
\end{eqnarray}
Though not explicit in the notation, the index set $I$ will always
have the above decomposition associated with it.  The variable $\iv$
shall denote the sequence $\iv=(|I_1|, |I_2|, \dots )$. We
denote by $l^2(I)$ the Hilbert space of square summable sequences
indexed by $I$ with inner product defined as $<\xv, \yv>= \sum_{i\in
I} x_i \bar{y_i} $.  We denote by $\delta _i$ the sequence whose $i$'th
entry is one and is zero otherwise; thus $\{ \delta _i \}_{ i \in I}$
is the canonical orthonormal basis for $l^2(I) $.

We shall let $\nu$ denote Lebesgue measure on $[0,1]$.

Let $\lfloor x \rfloor$ denote the greatest integer less than or equal to $x$.

Equality of two functions $f=g$ that have the same domain shall mean
that the two functions agree for every point in the domain.

Given two sequences $\xv= (x_1, x_2, \dots )$, $\yv=(y_1, y_2, \dots)$
and a scalar $c$, $\xv +\yv$ shall denote the sequence $(x_1 +y_1, x_2
+y_2, \dots )$, $c\xv $ shall denote the sequence $(cx_1, cx_2, cx_3,
\dots )$, $\frac{\xv}{\yv}$ shall denote the sequence
$(\frac{x_1}{y_1}, \frac{x_2}{y_2}, \dots)$, and $\lfloor \xv \rfloor$
shall denote the sequence $(\lfloor x_1 \rfloor, \lfloor x_2 \rfloor,
\dots )$.

$H$ shall denote a Hilbert space.  For a subset
$S\subset H$, $\span{S}$ shall denote the closure of the linear
subspace of $H$ spanned by the elements of $S$.  Given $h\in H$,
$||h||=(\ip{h}{h})^{\frac{1}{2}}$ shall denote the Hilbert space norm
of $h$.  Given $A:H \rightarrow H$, a bounded linear operator, 
$||A||=sup_{h \in H,\norm{h}=1} |\ip{Ah}{h}|$ shall be the operator
norm of $A$.

Appendix \ref{a:ultr} contains a summary of some basic notation and
properties of ultrafilters.

Finally, we remark that occassionally, when a result is
straightforward to verify, we will state it without providing a proof.
\subsection{Frames}

We use standard notations for frames as found in the texts of
Gr\"ochenig \cite{gr01}, or Daubechies
\cite{da92}; see also the research-tutorials \cite{hewa89} or
\cite{Cas00} for background on frames and Riesz bases.

We shall use the following  particular notation.

The definition of a {\it frame} is given in (\ref{framedef}).  A
sequence $\fc = \set{f_i}_{i \in I}$ that is a frame for $\span{\fc}$
which might not be all of $H$ shall be called a \emph{frame sequence}.

A frame is \emph{finite} if the size of the index set $I$ is finite
and \emph{infinite} if the size of the index set $I$ is
infinite. A frame is said to be {\em tight} if we can choose equal
frame bounds $A=B$. When $A=B=1$, the frame is called a {\em Parseval
frame}. We denote by $\fci$ the set of all frame sequences indexed by
$I$.

In the case of a frame or a frame sequence $\fc$, the \emph{frame
operator} $S$, defined by $Sf = \sum_{i \in I} \ip{f}{f_i} \, f_i$ is
a bounded, positive, and invertible mapping of $\span{\fc}$ onto
itself.  The \emph{Gram operator } $G$ in $l^2(I)$ is defined to be:
\begin{equation}
G: l^2(I) \rightarrow l^2(I), \ \ \ \{G(\{c_j\}_{j\in I})\}_i = 
\sum_{j\in I}<f_i, f_j> c_j.
\end{equation}

The following terminology is standardly applied to frames, however it
applies equally well to frame sequences; rather than introduce
additional notation, we shall associate to a frame or frame sequence
$\fc$:
\begin{itemize}
\item[-] the {\em canonical} (or \emph{standard}) \emph{dual
frame} $\tFc = \set{\tf_i}_{i \in I}$ where $\tf_i = S^{-1} f_i$.
\item[-] the \emph{associated Parseval frame} $\{
S^{-\frac{1}{2}}f_i \}_{i \in I}$ which has the property that it is
equal to its canonical dual frame and has upper and lower frame
bounds equal to 1.
\end{itemize}

The \emph{associated Gram projection} to a frame or frame sequence
$\fc$ will be the orthogonal projection in $l^2(I)$ onto the range of
the Gram operator $G$.  Equivalently, this is the Gram operator of the
associated Parseval frame.

A frame is a basis if and only if it is a Riesz basis, i.e., it is the
image of an orthonormal basis for $H$ under a continuous, invertible
mapping of $H$ onto itself. A Riesz sequence shall refer to a sequence
that is a Riesz basis for its closed linear span.

For two frames $\fc$ and $\gc$, the {\it superset} $\fc \oplus \gc$
shall denote the set $\{ f_i \oplus g_i\}_{i\in I}$.  Appendix
\ref{a:supersets} contains some basic notation and results pertaining
to supersets.

Note the upper bound inequality in (\ref{framedef}) is equivalent to
 $\norm{\sum_i c_i f_i}^2 \le B \, \sum_i |c_i|^2$ for any
 $(c_i)_i \in \ell^2(I)$.

\subsection{The sequences $a(\fc)$ and $b(\fc)$ associated to a frame.}

In this paper, frames will be compared using the data $\{f_i,
\tilde{f}_i\}_{i\in I}$.  Specifically, for each frame $\fc \in \fci$,
the sequence
\begin{equation}\label{eq:a}
a(\fc)= \{a_n(\fc)\}_{n \in \N}, \ a_n(\fc)= \frac{1}{|I_n|}
\sum_{i\in I_n} \ip{f_i}{\tilde{f}_i},
\end{equation}
shall play a central role.  The related ``unnormalized'' sequence
\begin{equation}
\label{eq:b}
b(\fc)= \{b_n(\fc)\}_{n \in \N}, \ b_n(\fc)= \sum_{i\in I_n}
\ip{f_i}{\tilde{f}_i}=|I_n|a_n(\fc),
\end{equation}
shall be used frequently.

\section{A new notion of frame equivalence} \label{s:equiv}

In this section we define the equivalence and partial ordering of
frames.  These concepts will only depend on the sequences $b(\fc)$ (
or equivalently $a(\fc)$).  The ideas and proofs about this
equivalence are more naturally viewed as properties of sequences.
Consequenctly we begin by defining a class of sequences, called {\it
frame compatible sequences} and showing that all sequences $b(\fc)$
arising from frames are frame compatible and that all frame compatible
sequences are "close" to $b(\fc)$ for some frame $\fc$ (Theorem
\ref{t:3.10}).  We then define an equivalence and partial order on
frame compatible sequences (Definition \ref{d:simm}) which naturally
pulls back to an equivalence and partial order of frames (Definition
\ref{d:frameequiv}).  We compare this equivalence to the well studied
standard equivalence.  Section \ref{s:3.3} shows the advantages of the
new equivalence.  Finally, in section \ref{s:3.4} we establish the
frame-sequence correspondence which relates the addition of sequences
to the superset operation $\oplus$ of certain frames (Theorem
\ref{t:frame-sequence}).  This correspondence will repeatedly be used
later in proofs about frame measure functions.

\subsection{Frame compatible sequences, equivalence and partial order}
\begin{Definition}  \label{d:framecompat} A sequence of nonnegative real
numbers $ \xv =(x_1, x_2, \dots )$ will be called {\it frame compatible} if 
\begin{enumerate}
\item $0 \leq x_1 \leq |I_1|$,
\item $0 \leq x_i- x_{i-1} \leq |I_i \backslash I_{i-1}|$ for all $i\geq2$.
\end{enumerate}
We shall denote by $\cs$ the set of all frame compatible sequences.
\end{Definition}
\begin{Remark}  Note that if $\xv$ is frame compatible then so is 
$\lfloor \xv \rfloor $. \end{Remark}
\begin{Definition}
A frame will be called {\it \fperp} if all nonzero elements of it 
are distinct elements of an orthonormal set.
\end{Definition}

\begin{Theorem} \label{t:3.10} \ \\
\begin{enumerate} 
\item Given a frame $\fc$, the sequence $b(\fc)$ is frame compatible.
\item For any frame compatible sequence $\xv$, there exists a \fperp 
frame denoted by $\gc ^{\xv}$ with $b( \gc ^{\xv}) = \lfloor \xv \rfloor$.

\end{enumerate}
\end{Theorem}
\pf Statement 1. follows simply from $0 \leq \ip{f_i}{\tilde{f}_i} \leq 1$. 

To prove $2.$ choose $S_1 \subset I_1$ such that $|S_1|=\lfloor x_1
\rfloor$.  Choose $S_i \subset I_i\backslash I_{i-1}$, $i \geq 2$ such
that $|S_i|=\lfloor x_i\rfloor - \lfloor x_{i-1} \rfloor$ (this can be
done precisely because $\xv$ is frame compatible).  Set $S = \cup_{i}
S_i \subset I$ and let $B$ be an arbitrary countable orthonormal
set.  Define a frame $\gc ^{\xv} = \{ g^{\xv}_i \}_{i \in I}$ such
that $\{g^{\xv}_i \}_{i \in S}$ are distinct elements of $B$ and
$g^{\xv}_i=0$ for $i \in I \backslash S$.  The frame $\gc^{\xv}$ is
normalized and tight since it is the union of distinct orthonormal
elements and zeroes.  It follows therefore that $<g^{\xv}_i,
\tilde{g^{\xv}_i}> = ||g^{\xv}_i||^2$ which is $1$ for $i \in S$ and
$0$ otherwise.  Thus $b_i(\gc^{\xv})= \sum_{j=1}^{i} |S_j|=\lfloor x_i
\rfloor$.  \qed

The following defines an important equivalence relation and partial
order on the set of frame compatible sequences.  We combine these
definitions with the map $b$ to produce the central object of this
paper: an equivalence and partial order on the set of frames $\fci$.

\begin{Definition} [Sequence equivalence and partial ordering]
\label{d:simm}
\ 

\begin{enumerate}
\item Given two sequences $\xv, \yv$ with complex entries we say
$\xv \simm \yv$ if $\lim_{n\rightarrow\infty} \frac{1}{|I_n|}(x_n - y_n)=0$.
\item Given two sequences $\xv$, $\yv$ with non-negative real entries
we say $\yv \leqq \xv$ if $\liminf_{n\rightarrow\infty}
 \frac{1}{|I_n|}(x_n - y_n) \geq 0$.
\end{enumerate}
\begin{Remark}  
For the moment, the equivalence relation and partial order will be
applied to frame compatible sequences.  However, later we shall be
considering this relation on a larger collection of sequences.
\end{Remark}
\end{Definition}
\begin{Definition}[Ultrafilter frame equivalence and partial ordering] 
\label{d:frameequiv}\mbox{}

\begin{enumerate}
\item We shall say two frames $\fc, \gc \in \fci$ are ultrafilter
equivalent, denoted $\fc \simm \gc $, if $b(\fc) \simm b(\gc)$.
\item For two frames $\fc$, $\gc \in \fci$ we say 
$\fc \leqq \gc$ if $b(\fc) \leqq b(\gc)$.\
\end{enumerate}
\end{Definition}

The next two subsections provide some motivation for this definition.
We begin by reviewing the standard notion of equivalence and then
discuss some advantages of the ultrafilter equivalence.

\subsection{The standard equivalence of frames.}

A different notion of equivalence of frames that has been studied
quite extensively is as follows (see \cite{ba99-2,hala00,allatawe04})
\begin{Definition}
 Given two frames $\fc =\{f_i\}_{i\in I}  \subset H_1$, 
$\gc =\{g_i\}_{i\in I}\subset H_2$, we say $\fc \sim \gc$ if 
there is a bounded invertible operator $S:H_1\rightarrow H_2$ 
such that $Sf_i=g_i$ for every $i\in I$.
\end{Definition}

It is easy to verify that $\sim$ is an equivalence relation (namely it
is reflexive, symmetric and transitive). Moreover, it admits the
following geometric interpretation that says that two frames are
$\sim$ equivalent if and only if the ranges of their Gram operators
are the same.
\begin{Theorem}[\cite{dailarson98,ba99-2}] \label{t.e1}
Consider $\fc ,\gc\ in\fci$ and let
$P,Q$ be their associated Gram projections.
Then $\fc \sim\gc$ if and only if $P=Q$.
\end{Theorem}
It is simple to verify that the equivalence in the $\sim$ relation
implies equivalence in the $\simm$ relation:
\begin{Proposition}\label{psim}
Given two frames $\fc$, $\gc \in \fci$, $\fc \sim \gc$ implies $\fc \simm \gc$.
\end{Proposition}

The $\sim$ equivalence relation is a very strong notion of
equivalence. For instance, in the following examples, the closely
related frames $\fc$ and $\gc$ are {\it not} $\sim$ equivalent.

\begin{Example} \label{e1}
 
Let the elements of $\gc$ differ from those in $\fc$ by scalars of
modulus one, i.e. $\gc=\{g_j=e^{i\phi_j} f_j: j\in I \}$.  In most
cases, these frames are not $\sim$ equivalent (unless $\fc$ was a
Riesz basis for its span).  In fact, this is true even when we require
that $e^{i \phi_j} \in \{-1,1 \}$.
\end{Example}
\begin{Example} \label{e2}
Let the elements of $\gc$ be a finite permutation of those in $\fc$,
i.e. let $\pi:I\rightarrow I$ be a finite permutation and set
$\gc=\{g_i=f_{\pi(i)};i\in I\}$.  In almost all cases $\fc$ and $\gc$
are not $\sim$ equivalent.
\end{Example}

\subsection{The advantages of the ultrafilter equivalence $\simm$.} 
\label{s:3.3}

The following proposition is strightforward and shows that unlike the
$\sim$ equivalence, the $\simm$ equivalence identifies the frames in
examples \ref{e1} and \ref{e2} as equivalent.
\begin{Proposition} \ 

\begin{enumerate}
\item If $\gc=\{g_j=e^{i\phi_j}f_j : j\in
I\}$, then  $\gc \simm \fc$.
\item If  $\gc=\{g_i=f_{\pi(i)} : i\in I\}$ for a finite permutation 
$\pi:I\rightarrow I$ then $\gc \simm \fc$.
\end{enumerate}
\end{Proposition}
The $\simm$ equivalence of frames holds for a much larger class of
permutations:

\begin{Proposition}
Let $\pi$ be a permutation (not necessarily finite) with the property that 
\[ \lim_{n \rightarrow \infty} \frac{|I_n \cap \pi (I_n)|}{|I_n|}=1.\]
If $\gc=\{ g_i = f_{\pi(i)} : \ i\in I\}$, then $\gc \simm \fc$.
\end{Proposition}

\pf Let $J_n= I_n \cap \pi (I_n)$, thus the sets $\{ f_j : j \in
J_n\}$ and $\{g_j : j \in J_n \}$ are identical.  The result follows
from the fact that:
\[ | a_n(\fc) - a_n(\gc)|\leq | a_n( \fc) - \frac{1}{|I_n|} \sum_{j \in J_n} \ip{f_j}{ \tilde{f}_j} |   +    |\frac{1}{|I_n|} \sum_{j \in J_n} \ip{g_j}{ \tilde{g}_j} - a_n(\gc)| \] \[\leq 2( 1- \frac{|J_n|}{|I_n|}), \]  the last inequality following from the fact that $\ip{f_j}{\tilde{f}_j},\ip{g_j}{ \tilde{g}_j}  \leq 1$. \qed

\bigskip
 At the heart of the ultrafilter equivalence is the sequence $a(\fc)$
(or equivalently $b(\fc)$).  Here we give an interpretation of
$a(\fc)$ from a stochastic signal analysis perspective.  This
interpretation further justifies the ultrafilter equivalence $\simm$.

\begin{figure}
\begin{center}
\mbox{\epsfig{figure=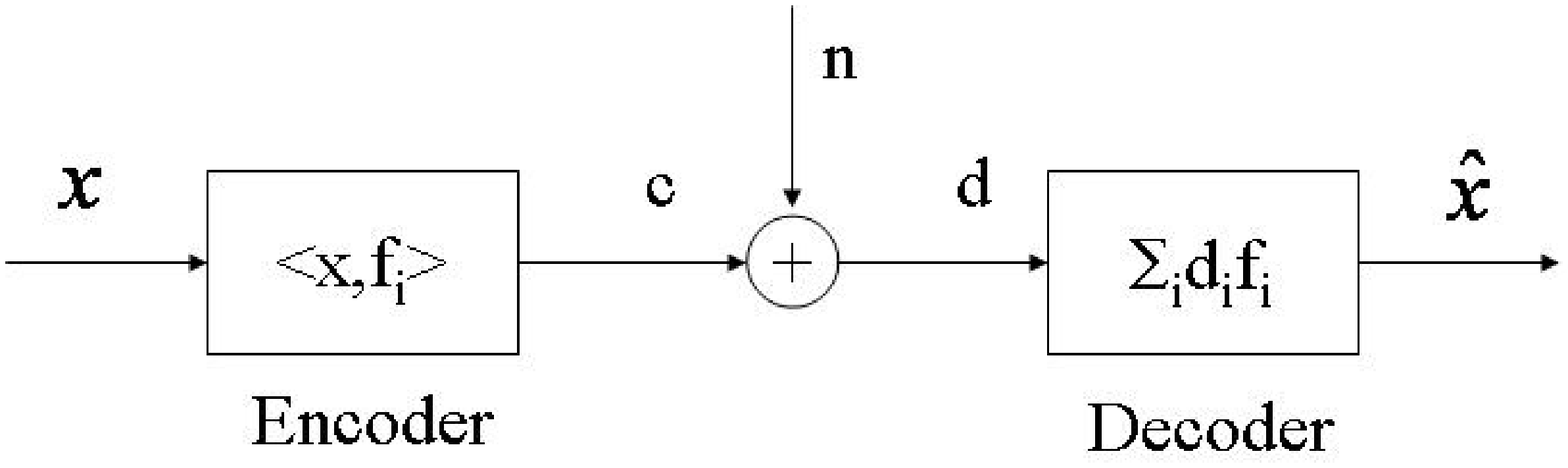,width=10cm,height=5cm,
bbllx=72,bblly=328, bburx=768,bbury=545}}
\end{center}
\caption{The Transmission Encoding-Decoding Scheme used to suggest the
importance of averages (\ref{e.e4}). \label{fig1}}
\end{figure}

We shall consider a Parseval frame $\fc\in\fci$.  Since every
frame is $\sim$ equivalent, (and thus $\simm$ equivalent by
Proposition \ref{psim}) to its associated Parseval frame, the behavior
of both equivalence relations is captured on the set of Parseval
frames.  Suppose the span $H$ of $\fc$ models a class of signals we
are interested in transmitting using an encoding and decoding scheme
based on $\fc$ as in Figure \ref{fig1}.

More specifically, a ``signal'', that is a vector $x\in H$, is
``encoded'' through the sequence of coefficients ${\bf c
}=\{\ip{x}{f_i} {\}}_{i\in I}$ given by the analysis operator
$T:H\rightarrow \ltwoi$. These coefficients are sent through a
communication channel to a receiver and there they are ``decoded''
using a linear reconstruction scheme $\hat{x}=\sum_{i\in I} d_if_i$
furnished by the reconstruction operator $T^{*}$. It is common to
consider what happens if the transmitted coefficients ${\bf
c}=(c_i)_{i\in I}$ are perturbed by some (channel) noise. In this
case, the received coefficients ${\bf d}=(d_i)_{i\in I}$ are not the
same as the transmitted coefficients ${\bf c}$. We shall assume the
system behaves as an {\em additive white noise channel} model, meaning
the transmitted coefficients are perturbed additively by unit variance
white noise.  Thus we can write
\begin{eqnarray}
d_i & = & c_i+n_i \label{e.e1} \\
\E[n_i] & = & 0 \label{e.e2} \\
\E[n_i\overline{n_j}] & = & \delta_{i,j} \label{e.e3}
\end{eqnarray}
where $\E$ is the expectation operator and $n_i$ represents the
independent noise component at the $i$'th coefficient. The
reconstructed signal $\hat{x}$ has two components, one due to the
transmitted coefficients $\sum_i c_if_i=x$ and the other due to the
noise $\eps=\sum_in_if_i$. We analyse the noise component.  Since its
variance is infinite in general, we consider the case that only
finitely many coefficients are transmitted, say a finite subset
$I_n\subset I$. Then the {\em average variance per coefficient of the
noise-due-error} is defined by:
\begin{equation}
a_n' = \frac{\E[|\eps _n|^2]}{|I_n|}
\end{equation}
where
\begin{equation}
\eps_n = \sum_{i\in I_n} n_if_i
\end{equation}
Using the assumptions (\ref{e.e2}),(\ref{e.e3}) we obtain
\begin{equation}
\label{e.e4}
a_n' = \frac{1}{|I_n|} \sum_{i\in I_n} \norm{f_i}^2
\end{equation}
which is exactly the quantity $a_n(\fc)$ used to define the
ultrafilter frame equivalence.  Since $\norm{f_i}\leq 1$ it follows
$a_n'\leq 1$. For an orthonormal basis the average noise-due-error
variance per coefficient would have been $1$ for all $n$ (since
$\norm{f_i}^2 =1$ for all $i$). Hence $a_n'=a_n(\fc)$ gives a measure
of how much the channel noise variance is reduced when a frame is used
instead of an orthonormal basis. In channel encoding theory, the noise
reduction phenomenon described before is attributed to the redundancy
a frame has compared to an orthonormal basis (see for instance
\cite{da90}). Hence, any measure of redundancy has to be connected to
the averages $a_n'=a_n(\fc)$ from (\ref{e.e4}).

It follows that two frames that are ultrafilter frame equivalent have
the same noise-due-error limiting behavior and if $\fc \leqq \gc$ then
$\fc$ has better noise-due-error limiting behavior.  The ultrafilter
frame measure function, which we introduce in section \ref{s:4.1}, is
defined using the limiting behavior of $a(\fc)$ to give an important
quantitative measure of frames.

\subsection{The frame sequence correspondence.}\label{s:3.4}

The following theorem describes the correspondence between frames and
frame sequences and shows that addition of frame sequences can be
realized by the superset operation ($\oplus$) of certain frames.

\begin{Theorem}[Frame-sequence correspondence] \label{t:frame-sequence}

\begin{enumerate}
\item
For every frame $\fc$ there exists a \fperp  frame $\gc$ with
$\lfloor b(\fc) \rfloor = b(\gc)$ and thus $\fc \simm \gc$.
\item Given frame compatible sequences $\xv^1, \dots , \xv^k$, and
$\zv=\sum_{i=1}^k \xv^i$, there exist frames $\fc ^{x^1}, \dots
\fc^{x^k}, \fc^z$ such that
\begin{enumerate}
\item $\fc^z= \oplus_{i=1}^k \fc^{x^i}$,
\item $b(\fc^{x^i})\simm \xv ^i$ for all $1\leq i \leq k$, and $b(\fc^z)\simm \zv$.
\end{enumerate}

\end{enumerate}
\end{Theorem}

{\it Proof of 1.:} Given $\fc$, the existence of $\gc$ is given by
Theorem \ref{t:3.10}.  It remains to show that the sequences $b(\fc)$
and $\lfloor b(\fc) \rfloor $ are $\simm$ compatible which follows
from
\[ \frac{b_n(\fc)-\lfloor b(\fc) \rfloor_n}{|I_n|} \leq \frac{1}{|I_n|} . \]

{\it Proof of 2.
:} 

We present the proof only for the case $k=2$. the general case follows
along same lines. 
We simplify the notation to $\xv^1=\xv$ and $\xv^2=\yv$.
Let $\bv$ be the sequence defined by: $\bv =
\lfloor \zv\rfloor - \lfloor\xv\rfloor -\lfloor\yv\rfloor$; notice
that $b_i\in\{0,1\}$.  Define $\tilde{\xv}=(\tilde{x_i})_i$ recursively
as:
\[ \tilde{x}_1= \lfloor x_1 \rfloor , \ \ \ \tilde{x_i}= min( \tilde{x}_{i-1} + \lfloor z_i \rfloor - \lfloor z_{i-1} \rfloor, \lfloor x_i \rfloor). \]
Using the fact that $\zv$ is frame compatible, it is straightforward
to verify that $\tilde{\xv}$ is frame compatible.  By definition,
$\tilde{x_i} \leq \lfloor x_i \rfloor $, we now show that $\tilde{x_i}
\geq \lfloor x_i \rfloor -1 + b_i$.  Suppose this is not the case then
let $j$ be the smallest index for which $\tilde{x}_{j}< \lfloor x_{j}
\rfloor -1 +b_j$.  Thus
\begin{eqnarray}
\tilde{x}_{j}&=& \tilde{x}_{j} +\lfloor z_j \rfloor - \lfloor z_{j-1} \rfloor \\
&=&\tilde{x}_{j-1}+ \lfloor x_{j} \rfloor -\lfloor x_{j-1} \rfloor +\lfloor y_{j} \rfloor - \lfloor y_{j-1} \rfloor + b_{j} - b_{j-1} \\
& \geq &\lfloor x_{j-1} \rfloor -1 +b_{j-1} +\lfloor x_{j} \rfloor -\lfloor x_{j-1} \rfloor +\lfloor y_{j} \rfloor - \lfloor y_{j-1} \rfloor + b_{j} - b_{j-1} \\
&=&\lfloor x_{j} \rfloor  -1 +\lfloor y_{j} \rfloor - \lfloor y_{j-1} \rfloor + b_{j} \geq \lfloor x_{j} \rfloor -1 +b_j
\end{eqnarray}
which contradicts the assumption on j.
Thus $\lfloor x_i \rfloor -1 \leq \tilde{x}_i \leq \lfloor x_i \rfloor
$ for all $i$, and hence $\tilde{\xv} \simm \lfloor \xv \rfloor \simm \xv$.  
Define
$\tilde{\yv} = \lfloor \zv \rfloor - \tilde{\xv}$.  It is
straightforward to verify from the definition of $\tilde{\xv}$ that
$\tilde{\yv}$ is frame compatible and since $\tilde{\xv} \simm \xv$,
we can conclude $\tilde{\yv} \simm \yv$.

By Theorem \ref{t:3.10}, since $\lfloor \zv \rfloor $ is frame
compatible, we can find a perpendicular-normal frame $\fc
^z=\{f_i^z\}_{i\in I} $ with $b(\fc ^z) = \lfloor \zv \rfloor $.
Define $T\subset I$ to be the subset of $I$ for which $f_i^z \not= 0$,
i.e. $T=\{ i \in I: f_i^z \not=0 \}$.  Write $T=T_1 \cup T_2$ such
that $T_1$ and $T_2$ are disjoint and $|T_1\cap I_i|=\tilde{x}_i$,
$|T_2 \cap I_i| = \tilde{y} _i$ ; this can be done since
$\tilde{\xv}$, $\tilde{\yv}$ and $\lfloor \zv \rfloor =\tilde{\xv} +
\tilde{\yv}$ are frame compatible.  Define $\fc ^x=\{ f^x_i\}_{i \in
I}$, $\fc^y = \{f^z_i\}_{i\in I}$ as follows: $f_i^x= f_i^z$ for $i\in
T_1$, $f_i^x =0$ otherwise, $f_i^y= f_i^z$ for $i\in T_2$, $f_i^y =0$
otherwise.  We have $\fc^x \oplus \fc ^y= \fc ^z$ (since the elements
of $\fc^z$ are orthogonal), and by construction $b(\fc^x)= \tilde{\xv}
\simm \xv$, $b(\fc^y)=\tilde{\yv} \simm \yv$ and $b(\fc^z)= \lfloor
\zv \rfloor \simm \zv$.

The proof of 3. follows along the same lines. \qed

\section{A measure of frames.} \label{s:ult}

In this section we introduce our main tool for a quantitative
comparison of frames: the ultrafilter frame measure function.  We give
its definition in section \ref{s:4.1}, and then we examine its
connection with the notion of excess in \ref{s:4.2}.  Appendix
\ref{a:ultr} gives a brief description of ultrafilters. Here we shall
denote by $N^*$ the set of free ultrafilters and for $p\in N^*$ and
$\xv =(x_1, x_2 , \dots )$ a sequence, the limit of $\xv$ along $p$
shall be denoted by $\plim\ \xv$. Finally $\cstarn$ shall denote the
set of continuous functions on $N^*$.

\subsection{The ultrafilter frame measure function}\label{s:4.1}
We shall now use ultrafilters to give a new measure for frames.

\begin{Definition}
Fix $(I_n)_{n\geq 0}$ as in Section \ref{s:2.1}. The {\em 
ultrafilter frame measure function} will be the map 
\begin{equation}
\label{eq.m3.1}
\mu:\fci\rightarrow\cstarn ~~;~~\mu(\fc)(p)=\plim \ a(\fc)=\plim \frac{1}{|I_n|}\sum_{i\in I_n}\ip{f_i}{\tilde{f}_i}  ,~~\forall p\in\N^*.
\end{equation}
\end{Definition}

\begin{Theorem}\label{t:4.2}

The ultrafilter frame measure function has the following properties:

\begin{enumerate}
\item $\mu(\fc _1) = \mu(\fc _2)$ if and only if $\fc _1 \simm \fc _2$.
\item $\mu(\fc _1)(p) \leq \mu(\fc _2)(p)$ for all $p \in \N^*$ if and only if $\fc _1 \leqq \fc_2$.

\item If $\fc$ is a Riesz basis for its span then $\mu(\fc)=1$.
\item If $\fc_1, \fc_2 \in \fci$ are such that $(\fc_1, \fc_2)$ are
orthogonal in the sense of supersets then $\mu(\fc_1 \oplus \fc_2)=
\mu(\fc_1) + \mu(\fc_2)$. (See Appendix \ref{a:supersets} for
definitions involving supersets.)
\end{enumerate}
\end{Theorem}
\pf

1. The statements
\begin{itemize}
\item[a)] $\mu(\fc _1) = \mu(\fc _2)$,
\item[b)] $\plim (a(\fc_1))= \plim(a(\fc_2))$ for all free 
ultrafilters $p$,
\item[c)] $\plim (a(\fc_1) -a(\fc_2))=0$ for all free ultrafilters $p$,
\item[d)] the sequence $a(\fc_1) -a(\fc_2)$ has a single accumulation point at $0$,
\item[e)] $\lim_{n\rightarrow \infty} (a(\fc_1)-a(\fc_2))=0$,
\item[f)] $\fc_1 \simm \fc_2$,
\end{itemize}
are all equivalent: $a) \Leftrightarrow b)$ and $e) \Leftrightarrow
f)$ follow from the definitions of $\mu$ and $\simm$, $b)
\Leftrightarrow c)$ follows from statement 2. of Proposition
\ref{p.m2.2}, $c) \Leftrightarrow d)$ is due to statement 3. of
Proposition \ref{p.m2.2}, $d) \Leftrightarrow e)$ follows from the
fact that $0\leq a_n(\fc_1), a_n(\fc_2) \leq 1$.

2. The proof is very similar to 1.; we omit the details.

3. If $\fc$ is a Riesz basis for its span,
$\ip{f_i}{\tilde{f_i}} =1$ for all $i\in I$.  Thus $a_n(\fc)=1$ for
all $n\in \N$ and so since $\lim _{n\rightarrow \infty}a_n(\fc)=1$,
statement 3. of Proposition \ref{p.m2.2} implies $\mu(\fc)=1$.

4. Since $\fc^1$ and $\fc^2$ are orthogonal in the
sense of supersets, the canonical dual frame of $\fc^1 \oplus \fc^2$
is $\{ \tilde{f_i^1} \oplus \tilde{f_i^2} \}_{i\in I}$, the direct sum
of the canonical duals for $\fc^1 $ and $\fc ^2$.  Since $\ip {f_i^1 \oplus f_i^2}{\tilde{f_i^1} \oplus \tilde{f_i^2}} =
\ip {f_i^1}{\tilde{f_i^1}} + \ip{f_i^2}{ \tilde{f_i^2}}$, we have
$a_n(\fc^1 \oplus \fc^2)= a_n(\fc^1) + a_n(\fc^2)$ and the result
follows. \qed

\subsection{The ultrafilter frame measure function and the excess of frames} \label{s:4.2}

The ultrafilter frame measure function gives information about the
excess of a frame -- a notion defined in \cite{bacahela03}.  We
begin by summarizing the relevant ideas and results of
\cite{bacahela03}.

The {\em excess} of a frame $\fc\in\fci$ with span $H$ is the supremum
over the cardinalities of all subsets $J\subset I$ so that
$\{f_i~:~i\in I\backslash J\}$ is complete in $H$. Since we consider
only countable sets $I$, the excess is either a finite number or
$\infty)$.  This supremum is always achieved \cite{bacahela03},
furthermore, for finite excess, $J$ can be always chosen so that
$\{f_i~:~i\in I\backslash J\}$ is also frame for $H$. However this
property no longer holds true in general for infinite excess. A
characterization of when this remains true was also given in
\cite{bacahela03}:
\begin{Theorem}[\cite{bacahela03}]\label{t.exc}
Let $\fc\in\fci$ be a frame for $H$ and $\tilde{\fc}$ its canonical dual. 
Then the following are equivalent:

a) There is an infinite subset $J\subset I$ such that $\{f_i~;~i\in
I\backslash J\}$ is frame for $H$;

b) There is an infinite subset $J'\subset I$ and $a<1$ so that
$\ip{f_i}{\tilde{f_i}} \leq a$ for all $i\in J'$.
\end{Theorem}

We now show that condition $b)$ is implied when the ultrafilter frame
measure function is not identically $1$.

\begin{Theorem}
\label{t.m.3}
Let $\fc\in\fci$ be a frame for $H$. If the ultrafilter frame measure
function $\mu(\fc)$ is not identically one, then there is an infinite
subset $J\subset I$ so that $\{f_i~;~i\in I\backslash J\}$ is frame
for $H$.
\end{Theorem}

\pf Since $\mu(\fc)$ takes on values in the interval
$[0,1]$, the hypothesis assumes that there exists some ultrafilter
$\accc$ such that $\mu(\fc)(\accc)<1$.  Thus we can find an infinite
set $J \in \accc$ and a constant $\epsilon >0$ such that
$a_j(\fc)<1-2\epsilon$ for all $j\in J$.  $a_j(\fc)$ is an average of
terms between 0 and 1 and thus it follows that at least
$\frac{\epsilon}{1-\epsilon}|I_j|$ of the terms
$\ip{f_i}{\tilde{f}_i}$, $i\in I_j$ are smaller than or equal to
$1-\epsilon$. Since $J$ is an infinite set, it follows that an
infinite number of the terms $\ip{f_i}{\tilde{f}_i}$ are bounded above
by $1-\epsilon$.  This establishes criterion $b)$ of \ref{t.exc} and
our result then follows. \qed

In subsequent papers \cite{bacahela03-1} and \cite{bacahela06-1} we analyzed
the excess problem for Gabor frames. There we showed that, if
the upper Beurling density is strictly larger than one then
there always exists an infinite subset that can be removed and leave
the remaining set frame. Furthermore, if the generating window belongs
to the modulation space $M^1$ and the lower Beurling density is
strictly larger than one, then one can find an infinite subset of
positive uniform Beurling density that can be removed and leave the
remaining set frame for $L^2$.

These results come as applications of the general theory we developed
in \cite{bacahela06}. There we analyzed the excess and overcompleteness for
a larger class of frames, namely those called localized frames. 
In that process
we obtained a completely new relation connecting the density of index
set to averages of the sequence $\{\ip{f_i}{\tilde{f}_i}\}$. We return
to this connection in Section \ref{s:8} in the context of Gabor
frames.  

Here we state one result from \cite{bacahela06} in our context.
To simplify notation, assume the index set $I$ is embedded in $\Z^d$,
that is $I\subset \Z^d$, for some integer $d$.
\begin{Definition}
A frame
$\fc\in\fc[I]$ is called {\em $l^1$-localized (with respect to its
canonical dual frame)} if there is a sequence $r\in\l^1(\Z^d)$ so that
$|\ip{f_i}{\tilde{f_j}}|\leq r(i-j)$.
\end{Definition}
For a subset $J\subset I\subset\Z^d$, we define its upper and lower 
{\em densities} as the following numbers:
\[ D^+(J) = \lim_{n\rightarrow \infty}\sup_{c\in\Z^d}\frac{|J\cap B_n(c)|
}{|B_n(c)|}
~~,~~D^-(J)=\lim_{n\rightarrow\infty}\inf_{c\in\Z^d}\frac{|J\cap B_n(c)|}{
|B_n(c)|}
\]
where $B_n(c)$ denotes the ball of radius $n$ centered at $c$ in $\Z^d$.
The set $J$ is said to have {\em uniform density} $D$ if $D^-(J)=D^+(J)=D$.
Now we restate Theorem 8 from \cite{bacahela06} using ultrafilter frame
measure function.
\begin{Theorem}
Assume $I\subset\Z^d$ for some integer $d$. Let $\fc\in\fc[I]$ be a 
$l^1$-localized frame for $H$. If $\mu(\fc)<1$ then there is an infinite
 subset $J\subset I$ of positive uniform density so that 
$\{f_i~;~i\in I\setminus J\}$ is frame for $H$.

Moreover, if $\mu(\fc)<\alpha<1$ then for each $0<\eps<1-\alpha$ the set
$J$ can be chosen as a subset of $\{i\in I ~;~\ip{f_i}{\tilde{f_i}}\leq
\alpha\}$ and the frame $\{f_i~;~i\in I\setminus J\}$ has a lower frame
bound $A(1-\eps-\alpha)$, where $A$ is the lower frame bound of $\fc$.
\end{Theorem}

\section{Sequence and frame measure functions} \label{s:gen}

The ultrafilter frame measure function provides a quantitative measure
for all frames indexed by the same set $I$.  In this section we
introduce the general notion of a frame measure function: a
quantitative measure of frames defined by some general properties
(Proposition \ref{framemeasurecharacterisation}).  We prove some
general facts about frame measure functions (Section \ref{s:5.1}) and
prove that the ultrafilter frame measure function has a lattice
structure (Section \ref{s:5.2}).  The natural way to view frame
measure functions is as linear maps on the sequences $a(\fc)$ via the
frame sequence correspondence (Theorem \ref{t:frame-sequence}).  For
this reason we present frame measure function via related maps on
sequences -- {\it sequence measure functions } (Definition
\ref{d:sequencemf}).  The technique of proving results about sequence
measure functions and "pulling the results back" to frame measure
functions will be used repeatedly through the rest of this work.

We begin by extending frame compatible sequences to a larger space 
of sequences.

\begin{Definition}
For the set $\cs$ of frame compatibe sequences, we let denote:
\begin{eqnarray}
\cs^+ & = & \{ c\xv ~:~\xv\in\cs~,~c\geq 0\} \\
\cse & = & \{ \xv^1-\xv^2 ~:~\xv^j\in\cs^+~{\rm for}~1\leq j\leq 2\}
\end{eqnarray}
\end{Definition}
 
\begin{Proposition}
\ 

\begin{enumerate}
\item The set $\cs$ of frame compatible sequences is convex.
\item If $0\leq c\leq 1$ and $\xv\in\cs$, then $c\xv\in\cs$.
\item $\cs^+$ is a positive cone, that is, for $c_1,c_2\geq 0$, and
 $\xv_1,\xv_2\in\cs^+$, we have $c_1\xv_1+c_2\xv_2\in\cs^+$.
\item $\cse$ is the real vector space spanned by $\cs$, that is for
any $c_1,c_2\in\R$, $\xv_1,\xv_2\in\cse$, we have $c_1\xv_1+c_2\xv_2\in\cse$.
\end{enumerate}
\end{Proposition}

\pf

Property 1. is a consequence of the fact that the constraints of
the definition of frame compatibility (Definition \ref{d:framecompat})
are convex.  Property 2.  follows from convexity of $\cs$, since both
$0$ and $\xv$ belong to $\cs$.  Property 3. follows from 1. and 2.
Finally property 4. follows from definition of $\cse$ and 3.  \qed

\begin{Theorem} \label{t:uniqueextension}

Given a linear function $m$ on the frame compatible sequences, there
exists a unique linear extension $\tm$ of $m$ to $\rs$.

\end{Theorem}

\pf Since $m$ is linear on $\cs$, it is clear that defining $\tm (c
\xv) \equiv cm(\xv)$ for $\xv$ frame compatible and $c \geq 0$
uniquely extends $m$ to $\cs^+$.  Linearity of $m$ on $\cs$ implies
linearity of $\tm$ on $\cs^+$ as follows: for $\xv, \ \yv \in \cs$,
$c,d >0$ we have
\[ \tm (c \xv + d \yv)= (c+d)m(\frac{c}{c+d} \xv + \frac{d}{c+d} 
\yv)=(c+d)(\frac{c}{c+d}m( \xv) + \frac{d}{c+d} m(\yv))
=c m(\xv) + d m(\yv),\]
since $\frac{c}{c+d}\xv, \ \frac{d}{c+d}\yv, \frac{c}{c+d}\xv + 
\frac{d}{c+d}\yv \in \cs$. If a linear extension to $\cse$ existed, 
it would have to be unique since $\xv \in \cse$ implies 
$\xv =\xv ^1 - \xv ^2 $ for some $\xv^{j} \in \cs^+, \   
1\leq j \leq 2$.  Hence by linearity we would have to have 

\begin{equation}\label{e:complexdef}
\tm(\xv) = \tm (\xv^1)- \tm (\xv ^2) .
\end{equation}

It remains to show that (\ref{e:complexdef}) is well defined.  Suppose
$\xv= \xv ^1 - \xv ^2 =\yv ^1 - \yv ^2 $ for $\xv^j, \yv^j \in \cs
^+$, $ 1\leq j \leq 2$.  Then $\xv ^1 + \yv ^2=\yv ^1 + \xv ^2 $.  By
the linearity of $\tm$ on $\cs^+$ we have $\tm (\xv ^1) + \tm (\yv ^2) =
\tm (\yv ^1) + \tm (\xv ^2)$.  Rearranging terms yields
\[\tm (\xv^1)- \tm (\xv ^2) = \tm (\yv^1)- \tm (\yv ^2)  , \] and thus  
(\ref{e:complexdef}) is well defined. \qed

\begin{Definition}

Let $W$ be a compact Hausdorff space; denote by $\cstarw$ the set of
real-valued continuous functions over $W$.
\end{Definition}

We now define the notions of a sequence and frame measure function.

\begin{Definition}\label{d:sequencemf}
 A sequence measure function $m: \cse \rightarrow \cstarw$ will be a
 function which satisfies
\begin{enumerate}
\item For $x$, $y \in \cse$, 
$m( \xv) = m( \yv)$ if and only if $\xv \simm \yv$,
\item For $x$, $y \in X^+$, $m( \xv) \leq m( \yv)$ if and only if $\xv
\leqq \yv$,
\item For $\iv =( |I_1|, |I_2|, |I_3|, \dots )$, $m(\iv)=1$,
\item $m$ is linear.
\end{enumerate}
\end{Definition}

\begin{Definition} 

A {\it frame measure function} will be a function
$m_{f}:\fci\rightarrow\cstarw$ which is the composition of the map $b
: \fci \rightarrow X$ and a sequence measure function $m$,
i.e. $m_f(\fc)=m(b(\fc))$, for all $\fc \in \fci$.
\end{Definition}
The ultrafilter frame measure function is a frame measure function as
we prove in Corollary \ref{cor-X}.

An equivalent description of a frame measure function is as follows:
\begin{Proposition} \label{framemeasurecharacterisation}
 A map $m_f:\fci\rightarrow\cstarw$ is a frame measure function if and
 only if it satisfies the following properties:
\begin{enumerate}
\item[A.] $m_f(\fc _1) = m_f(\fc _2)$ if and only if $\fc _1 \simm \fc _2$.
\item[B.] $m_f(\fc _1)(x) \leq m_f(\fc _2)(x)$ for all $x \in M$ if
and only if $\fc _1 \leqq \fc_2$.
\item[C.] If $\fc$ is a Riesz basis for its span then $m_f(\fc)=1$.
\item[D.] If $\fc_1, \fc_2 \in \fci$ are such that $(\fc_1, \fc_2)$
are orthogonal in the sense of supersets then $m_f(\fc_1 \oplus
\fc_2)= m_f(\fc_1) + m(\fc_2)$.
\end{enumerate}

\end{Proposition}
\pf Given a frame measure function $m_f=m\circ b$, properties A. and
 B. follow immediately from properties 1. and 2. of Definition
 \ref{d:sequencemf}.  Given a Riesz basis for its span $\fc$, we have
 $\ip{f_i}{ \tilde{f}_i}=1$ for all $i \in I$ and hence
 $b(\fc)=\{|I_1|, |I_2|, \dots \}$.  Thus property C. above follows
 from property 3. of Definition \ref{d:sequencemf}.  Finally, if $\fc
 _1$, $\fc _2$ are orthogonal in the sense of supersets we have $b(\fc
 _1 \oplus \fc _2)= b(\fc_1) + b(\fc _2)$ and the linearity (property
 4.) of $m$ implies property D. above.
 
 We are left to show that a map $m_f$ satisfying the above 4
 properties implies that the existence of a sequence measure function
 $m$ with $m_f= m \circ b$.  We first define $m$ on the frame
 compatible sequences from $m_f$ as follows.  Given $\xv \in \cs $, by
 Theorem \ref{t:3.10} there is a frame $\gc^{\xv}$ with
 $b(\gc^{\xv})=\lfloor \xv \rfloor$, we define $m(\xv)=m_f(\gc
 ^{\xv})$.  Now for any frame $\fc$, if we let $\xv=b(\fc)$, we have
 $\fc \simm \gc^{\xv}$ since $b(\fc) \simm \lfloor \xv \rfloor
 =b(\gc^{\xv})$.  Thus by condition A.,
 $m_f(\fc)=m_f(\gc^{\xv})=m(\xv)=m(b(\fc))$ and thus $m_f=m \circ b$.
 
 We now show this map $m$ is linear on the set of frame compatible
 sequences, i.e.
\begin{enumerate}
\item if $\xv$ and $c \xv$ are frame compatible then $cm(\xv)=m(c \xv)$,
\item if $\xv$, $\yv$ and $\xv + \yv$ are frame compatible then 
$m(\xv) + m(\yv)= m(\xv + \yv)$.
\end{enumerate}

For any $\frac{a}{b} < c$, $a,\ b \in \N$, set $\yv= \frac{1}{b} \xv
\in \cs$.  Applying part 3. of Theorem \ref{t:frame-sequence} to the
case $k=b$, $\xv^i=\yv$, $1\leq i \leq k$ yields $bm(\yv)=m(\xv)$.
Similarly $am(\yv)= m(\frac{a}{b} \xv)$; combining these conditions
yields $\frac{a}{b} m(\xv)= m(\frac{a}{b}\xv)$.  Since $\frac{a}{b}
\xv \leqq c\xv$, properties A. and B. imply $m(\frac{a}{b}\xv) \leq
m(c \xv)$.  Coupling this with the above two relations yields
$\frac{a}{b} m(\xv) \leq m(c\xv)$.  Applying this to a sequence of
rational $\frac{a}{b}$ that approach $c$ from below yield $c m(\xv)
\leq m(c\xv)$.  A similar argument can be made for any rational
fraction greater than or equal to $c$ and we conclude $cm(\xv) \leq
m(c\xv) \leq cm(\xv)$ and thus $cm(\xv)=m(c\xv)$.

Statement 2. above follows directly from property D. and part 2. of
Theorem \ref{t:frame-sequence}.

Thus $m$ is linear on the set of frame compatible sequences and by
Theorem \ref{t:uniqueextension} we can uniquely extend $m$ to a linear
map on $\rs$; we will call this extended map $m$ as well.  It remains
to show that $m$ satisfies properties $1.-3.$ of Definition
\ref{d:sequencemf}.  Property $3.$ follows from the fact that for an
orthonormal basis $\fc$, $m_f(\fc)=1$ and $b(\fc)=\{ |I_1|, |I_2|,
|I_3|, \dots \}$.  We now establish property $1.$ of Definition
\ref{d:sequencemf}. Given $\xv,\yv \in \cse$, write $\xv =\xv ^1 - \xv
^2 $, $\yv =\yv ^1 - \yv ^2 $, with $\xv^j, \yv^j \in X^+$ and
$\frac{1}{c}\xv^j$, $\frac{1}{c}\yv^j$ frame compatible sequences.  It
is straightforward to verify that

\[ m(\xv) = m (\yv) \Leftrightarrow \tm (\xv^1 +\yv^2)=m (\yv^1 +\xv^2),  \]

\[ \Leftrightarrow m(\frac{1}{2c}(\xv^1 +\yv^2))=m(\frac{1}{2c}(\xv^2 +
\yv^1))\]

\[ \Leftrightarrow \frac{1}{2c}(\xv^1 +\yv^2)\simm \frac{1}{2c}(\xv^2 +\yv^1)\]

\[ \Leftrightarrow \xv^1 +\yv^2\simm \xv^2 +\yv^1\]

\[\Leftrightarrow \xv \simm \yv \]

where the third double implication comes from property 1. of a frame
measure function and all other implications follow from the linearity
of $m$.  Finally we show property $2.$ of Definition
\ref{d:sequencemf}. Given $\xv,\yv \in X^+$, there exist a constant
$c$ such that $\frac{1}{c}\xv, \ \frac{1}{c}\yv $ are both frame
compatible.  It is then straightforward that

\[ m(\xv) \leq m(\yv) \Leftrightarrow m(\frac{1}{c}\xv) \leq m(\frac{1}{c}\yv)  \Leftrightarrow (\frac{1}{c}\xv) \leqq (\frac{1}{c}\yv) \]

\[ \Leftrightarrow \xv \leqq \yv. \qed \]

\begin{Remark}  
Condition $D.$ in Proposition \ref{framemeasurecharacterisation} can
be viewed as a linearity condition on supersets of certain pairs of
frames.  One might hope for more, namely that one could find a map
with conditions $A$, $B$ and $C$ with the added property that the map
was linear on supersets of all pairs of frames.  This turns out to be
too much to hope for as the following example shows:

\begin{Example} 

Let $H$ be a Hilbert space with orthonormal basis $\{e_i\}_{i\in \N}$,
let $I_n=\{1, \dots , n\}$.  Define $\fc=\{f_i\}_{i\in \N}$ and
$\gc=\{g_i\}_{i \in \N}$ as follows:

\[ f_i= \left\{ \begin{array}{rcl} e_i & , & i \mbox{ even} \\ 
0 & , & i \mbox{ odd.} \end{array}\right. \]
\[ g_i= \left\{ \begin{array}{rcl} 
\frac{1}{2}e_i + \frac{1}{2} e_{i^2 +1} & , & i \mbox{ even} \\  
\frac{1}{2} e_{\sqrt{i-1}} + \frac{1}{2} e_i & , &  \sqrt{i-1} \mbox{ even.} \\
0 & , & \mbox{ otherwise.} \end{array}\right. \]

Let $H_1= \span{\fc}$, $H_2= \span{\gc}$.  The following facts about
$\fc$ and $\gc$ can be verified:
\begin{itemize}
\item $\span{\{ f_i \oplus g_i\}}_{i\in \N}= H_1 \oplus H_2$,
\item $\fc \oplus \gc =\{f_i \oplus g_i\}_{i \in \N}$ is a frame for $H_1 \oplus H_2$ (this is verified using Theorem \ref{p:A2} or checking that $h_i$ below are the dual frame elements),
\item the canonical dual frame $\{ \tilde{h}_i\}$ is given by 
\[ \tilde{h}_i= \left\{ \begin{array}{rcl} e_i \oplus 0 & , & i \mbox{ even} \\
 -e_{\sqrt{i-1}} \oplus e_{\sqrt{i-1}} + e_i & , & \sqrt{i-1}  \mbox{ even} \\ 
0 & , & \mbox{ otherwise.} \end{array}\right. \]
\item $\mu(\fc)=\frac{1}{2}$, $\mu (\gc)= \frac{1}{4}$, $\mu(\fc \oplus \gc)= \frac{1}{2}$.
\end{itemize}

Thus $\mu$ is not additive in the sense of supersets in this case.

Though this shows no map of the above form can be linear on supersets
of all pairs of frames, a main result of Section \ref{s:9} shall be
that for index sets $I$ with a little added structure, frame measure
functions are linear on supersets of pairs of frames coming from a
large subset of all frames that includes Gabor frames.
\end{Example}
\end{Remark}

It is straightforward to verify using Proposition
\ref{framemeasurecharacterisation} that: 
\begin{Corollary}\label{cor-X}
The ultrafilter frame measure function is a frame measure function.
\end{Corollary}

We define the corresponding sequence measure function:
\begin{Definition}
The {\em ultrafilter sequence measure function} shall be the
sequence measure function corresponding to the ultrafilter frame
measure function, i.e. the map
\[ 
\mu: \ \cse \rightarrow \cstarn \ \ \ \mbox{ given by } \ \ \ 
\mu(\xv)(p)= \plim \frac{x_n}{|I_n|} . \]
\end{Definition}
We will use the same $\mu$ to denote both the ultrafilter
sequence measure function and the ultrafilter frame measure
function.

\subsection{General properties of sequence and frame measure functions.} 
\label{s:5.1}
\begin{Proposition} \label{p:sequenceprops}
Suppose $\xv \in \cse $ and $m$ is a sequence measure function, then
\begin{enumerate}
\item If $c=\lim_{i \rightarrow \infty} \frac{x_i}{|I_i|}$ exists then
$m(\xv)$ is the constant function of value $c$.
\item $\liminf \frac{x_i}{|I_i|} \leq m(\xv)(w) \leq \limsup
\frac{x_i}{|I_i|}$ for all $w\in W$.
\item There exist $v,w\in W$ (different for different $\xv$) such that
$m(\xv)(v)=\liminf \frac{x_i}{|I_i|}$, $m(\xv)(w)= \limsup
\frac{x_i}{|I_i|}$ .

\end{enumerate}

\end{Proposition}
{\it Proof of 1}: Recall $\iv= (|I_i|)_i$.  Set $\yv= c \iv$.  It
follows from Definition \ref{d:simm} that $\xv \simm \yv$ and so
$m(\xv)=m(\yv)=cm(\iv)=c\cdot 1$, the last two equalities following from
the linearity of $m$ and condition 3 of Definition \ref{d:sequencemf}.

{\it Proof of 2. and 3}: Let $l$ be the greatest number for which
$l\iv \leqq \xv$ and let $L$ be the smallest number for which $\xv
\leqq L\iv$.  From the definition of $\leqq$, it follows that $l=
\liminf \frac{\xv}{\iv}$ and $L= \limsup \frac{\xv}{\iv}$; result 2. then
follows since $m(c\iv)(w)=c$ for all $w\in W$.  Furthermore,
property 2. of the definition of a sequence measure function
(Definition \ref{d:sequencemf})) ensures that $l=\liminf_{w\in W}
m(\xv)(w)$ and $L=\limsup_{w\in W} m(\xv)(w)$.  The continuity of
$m$ and the compactness of $W$ ensures that there exist points $v,w
\in W$ for which $m(\xv)(w)$ achieves the lower and upper bounds,
i.e. $m(x)(v)=l$, $m(x)(w)= L$.  \qed

Proposition \ref{p:sequenceprops} "pulls back" via the map $b$ and the
frame sequence correspondence (Theorem \ref{t:frame-sequence}) to
the following statement about frame measure functions:
\begin{Proposition} \label{p:frameprops}
Suppose $\fc \in \fci $ and $m$ is a frame measure function, then
\begin{enumerate}
\item If $c=\lim_{i \rightarrow \infty} a_i(\fc)$ exists then $m(\fc)$
is the constant function of value $c$.
\item $\liminf a(\fc) \leq m(\fc)(w) \leq \limsup a(\fc)$ for all $w\in W$.
\item There exist $v,w\in W$ (different for different $\fc$) such that
$m(\fc)(v)=\liminf a(\fc)$, $m(\fc)(w)= \limsup a(\fc)$ .

\end{enumerate}
\end{Proposition}

\subsection{Sequences and Lattices} \label{s:5.2}

\begin{Proposition} \label{p:sequencechar}
A real valued sequence $\xv$ is in $\rs$ if and only if there exists a
constant $c$ such that $|x_1|\leq c|I_1|$ and $|x_{i}-x_{i-1}| \leq
c(|I_{i}| - |I_{i-1}|)$ for all $i\geq 2$.
\end{Proposition}
\pf If $\xv \in \rs$ then $\xv = \xv ^1 - \xv ^2$ with $\xv ^1, \xv ^2
\in \cs ^+$.  Thus there exists a constant $c$ such that
$\frac{2}{c}\xv^1$, $\frac{2}{c} \xv ^2 \in \cs$ and therefore, $|x
^k_1|\leq \frac{c}{2} |I_1|$ and $|x ^k _{i} - x^k_{i-1}| \leq
\frac{c}{2} (|I_i|- |I_{i-1}|)$ for $k=1,2$.  It follows that
$|x_1|\leq c |I_1|$ and $|x _i- x _{i-1}| \leq c(|I_i|-
|I_{i-1}|)$.

Given a sequence $\xv$ such that there is a constant $c$ for which
$|x_1| \leq c |I_1|$, $|x_i- x_{i-1}| \leq c(|I_i|- |I_{i-1}|)$,
$i\geq 2$, set $d_1=x_1$, $d_i= x_i - x_{i-1}$ for $i\geq 2$.
Inductively define $\xv ^1$, $\xv ^2$ as follows: \ \ $ x^1_1= \max
(d_1, 0), \ x^2_1 = \max (-d_1, 0), \ \ \ x^1_{i}= x^1_{i-1} + \max
(d_i, 0), \ \ \ x^2_{i}= x^2_i + \max(-d_i, 0).$ By construction $\xv
^1 - \xv ^2= \xv$.  In addition, $x^k_i- x^k_{i-1} \in \{0, |d_i|\}$
and thus $\xv ^1$, $\xv ^2 \in \cs ^+$. \qed

\begin{Definition} \label{d:wedgevee}
For real valued sequences $\xv$, $\yv$, define the sequences $\xv
\wedge \yv$ and $\xv \vee \yv$ as follows:
\[ (\xv \wedge \yv)_i = \min (x_i, y_i) \hspace{.7in} (\xv \vee \yv)_i = \max (x_i, y_i) \mbox{   \ \ \ \ for all $i\geq 1.$} \]
\end{Definition}
\begin{Remark} \label{r:veewedge}
It follows from the definitions that $(\xv \wedge \yv) \leqq \xv \leqq
(\xv \vee \yv)$ and $(\xv \wedge \yv )\leqq \yv \leqq (\xv \vee \yv) $.
\end{Remark}
\begin{Proposition} \label{p:lattice}
The sets $\cs, X^+, \cse$ are all closed under the binary operations
$\wedge$ and $\vee$. Consequently each set forms a lattice.
\end{Proposition}
\pf Given $\xv, \yv \in \cs $ and $i\geq 1$, without loss of
generality we can assume $(\xv \wedge \yv)_{i-1} =x_{i-1}$.  So
\[ (\xv \wedge \yv)_i- (\xv \wedge \yv)_{i-1} = \min (x_i, y_i) - x_{i-1} 
\leq x_i-x_{i-1} \leq |I_i|- |I_{i-1}|, \] 
and so $\xv \wedge \yv \in \cs$.  The result for $\cs ^+$ follows from 
the result for $\cs $ by noting that $c (\xv \wedge \yv)= c\xv \wedge c \yv$.

For $\xv, \yv \in \rs$ let $c$ be as in proposition
\ref{p:sequencechar} so that $|x_1|, \ |y_1| \leq c |I_1|$ and $|x_i-
x_{i-1}|, |y_i - y_{i-1}| \leq c ( |I_i|- |I_{i-1}|)$ for $i\geq 2$.
We now consider the two cases a) $(\xv \wedge \yv)_i \geq (\xv \wedge
\yv)_{i-1}$, b) $(\xv \wedge \yv)_i <(\xv \wedge \yv)_{i-1}$.  In case
a) we can assume $(\xv \wedge \yv)_{i-1} = x_{i-1}$ and again
\[ 0 \leq (\xv \wedge \yv)_i  - (\xv \wedge \yv)_{i-1}= \min (x_i, y_i) 
- x_{i-1} \leq x_i - x_{i-1} \leq c(|I_i|- |I_{i-1}|).\]  
In case b) we can assume $(\xv \wedge \yv)_i= x_i$ and thus
\[ 0 \geq (\xv \wedge \yv)_i -(\xv \wedge \yv)_{i-1}= 
x_i -\min (x_{i-1}, y_{i-1})  \geq x_i - x_{i-1} \geq -c(|I_i|- 
|I_{i-1}|) .\]  
These two cases establish that $\xv \wedge \yv$ satisfy the 
conditions of Propositon \ref{p:sequencechar} and thus $\xv 
\wedge \yv \in \rs$.

The corresponding result for $\xv \vee \yv$ can be proven in a similar
fashion. \qed

\begin{Proposition} \label{p:minmax}
The ultrafilter sequence measure function has the properties:
\begin{enumerate}
\item $\mu (\xv \wedge \yv) (p)= \min(\mu(x) (p), \mu(y)(p))$,
\item  $\mu (\xv \vee \yv) (p)= \max (\mu(x) (p), \mu(y)(p))$.
\end{enumerate}
\end{Proposition}

The lattice structure on sequences induces a lattice structure on frames:

\begin{Definition}
Given two frame $\fc$, $\gc$, $\fc \vee \gc$ will denote any frame
that has the property that $b(\fc \vee \gc)\simm b(\fc) \vee b(\gc)$.
Similarly, denote by $\fc \wedge \gc$ any frame that has the property
that $b(\fc \wedge \gc) \simm b(\fc) \wedge b(\gc)$.
\end{Definition}
\begin{Remark}
Theorem \ref{t:frame-sequence} and Proposition \ref{p:lattice}
guarantee the existence of the frames $\fc \wedge \gc$ and $\fc \vee
\gc$.
\end{Remark}

With this notation Proposition \ref{p:minmax} implies:
\begin{Proposition} \label{prop.5.20}
The ultrafilter frame measure function has the properties:
  
 \begin{enumerate}
\item $\mu (\fc \wedge \gc) (p)= \min(\mu(\fc) (p), \mu(\gc)(p))$,
\item  $\mu (\fc \vee \gc) (p)= \max (\mu(\fc) (p), \mu(\gc)(p))$.
\end{enumerate}
\end{Proposition}

\subsection{Universality of the ultrafilter sequence and frame 
measure function}
We now show that a copy of the ultrafilter sequence measure function
is embedded in any sequence measure function and consequently a copy
of the ultrafilter frame measure function is embedded in any frame
measure function.
\begin{Theorem}\label{t:canonembed}
Given a sequence measure function $m$, and an ultrafilter $\accc$,
there exists an element $w_{\accc} \in W$ such that $\mu (\xv)(\accc)=
m(\xv)(w_{\accc})$ for all $\xv \in \cse$.
\end{Theorem}

\pf Given an ultrafilter $\accc$, denote by $Y_{\accc}$ all
sequences for which the ultrafilter limit along $\accc$ is the
$\limsup$ of the sequence, i.e.

\[Y_{\accc}= \{ \yv \in \rs: \mu(\yv)(\accc)= \limsup \frac{\yv}{\iv} \}.
\]
  Set $W_{\accc}= \{ w \in W: m(\yv)(w)=\limsup \frac{\yv}{\iv}
\mbox{ for all } \yv \in Y_{\accc} \}$.  We will eventually show that
every point $w \in W_{\accc}$ satisfies $m(\xv)(w)= \mu(\xv)(\accc)$
for all $\xv \in \cse$.  We begin by showing that $W_{\accc}$ is
nonempty.

\begin{Lemma}
For all free ultrafilters $\accc$, $W_{\accc}$ is nonempty.
\end{Lemma}

\pf
Suppose $W_{\accc}= \emptyset$ for some $\accc$.  Thus for every point
$w\in W$ there exists a sequence $\yv^w \in Y_{\accc}$ such that
$m(\yv^w)(w) < \mu (\yv^w)(\accc)= \limsup \frac{\yv^w}{\iv}$. 
Since $m$ is
continuous we can find an open set $V_w$ around $w$ such that
$m(\yv^w)(v) \leq c_w < \mu (\yv^w)(\accc)$ for all $v \in V_w$.  Thus
$\cup_{w\in W}V_w$ is an open cover of $W$.  Since $W$ is compact we
can find $w_1 , \dots, w_n$ such that $\cup_{i=1}^n V_{w_i}=W$ and
therefore for all $w\in W$ there exists an $i(w) \in \{1,2, \dots,
n\}$ such that

$m(\yv^{w_{i(w)}})(w) \leq c_{w_{i(w)}} < \mu  (\yv^{w_{i(w)}})(\accc)$.  Setting $\zv= \sum_{i=1}^n \frac{1}{n} \yv^{w_i}$ we have

\[m(\zv)(w) = \sum_{i=1}^n \frac{1}{n} m(\yv^{w_i})(w )   \leq \sum_{i\not= i(w)}\frac{1}{n} m(\yv^{w_i})(w) + \frac{1}{n}c_{w_{i(w)}} < \frac{1}{n} \sum_{i=1}^n \mu  (\yv^{w_i})(p)= \mu  (\zv)(p)\] 

for all $w$.  This however contradicts Proposition
\ref{p:sequenceprops} since it shows that $m(\zv)$ cannot achieve
$\limsup \frac{z_n}{|I_n|}$ since it is strictly less than $\mu (\zv
)(p)$.\qed

The lemma established that $W_{\accc}$ is nonempty; we now show that
each $w \in W_{\accc}$ has the property that $m(\xv)(w)=\mu
(\xv)(\accc)$ for all $\xv \in \cse$.  Suppose this is not the case,
i.e. there is an $\xv$ such that $m(\xv)(w) \not= \mu (\xv)(\accc)$.
Assume first that $m(\xv)(w) < r < \mu (\xv)(\accc)$.  Set $\yv= \xv
\wedge r \iv $ ( see Definition \ref{d:wedgevee} ).  Remark
\ref{r:veewedge} then implies that $m(\yv)(w) \leq m(\xv)(w) < r$.  In
addition $\mu (\yv)(\accc)=r$ by Proposition \ref{prop.5.20}.  
However, since 
$$r=\plim \frac{\yv}{\iv} \leq \limsup (\frac{\yv}{\iv}) = 
\limsup (\min (\frac{x_n}{|I_n|}, r))\leq r$$
 we have
$\yv \in Y_{\accc}$ and thus by the definition of $W_{\accc}$ we must
have $m(\yv)(w)=r$, a contradiction.

The case $m(\xv)(w) > \mu (\xv) (\accc)$ reduces to the previous case
by noting that for $\xv' = \iv - \xv$ we have $m(\xv')(w)=1-
m(\xv)(w)< 1 - \mu(\xv)(\accc)=\mu(\xv')(\accc)$.  \qed

The following corollary follows from the frame-sequence correspondence
(Theorem \ref{t:frame-sequence}):

\begin{Corollary} \label{c:frameembed}
Given a frame measure function $m$, and an ultrafilter $\accc$, there
exists an element $w_{\accc} \in W$ such that $\mu (\fc)(\accc)=
m(\fc)(w_{\accc})$ for all $\fc \in \fci$.
\end{Corollary}

\section{Topological results} \label{s:top}
We now examine sequence and frame measure functions from a topological
point of view.
 
Corollary \ref{c:frameembed} says that a copy of the ultrafilter frame
measure function $\mu $ can be found inside any frame measure
function.  However, this is only an algebraic copy and nothing has
been shown about the topological compatibilities between the two
measure functions.  We partially address these issues in this section.
In \ref{subs:top1} we introduce some natural additional properties
(separable, irreducible, minimal) that a sequence or frame measure
function could have and we define a canonical minimal measure function
$\mu ^0$ related to $\mu$.  We also give a canonical construction for
turning an arbitrary sequence or frame measure function into a
separable one.  In \ref{subs:top2} we prove two important results:
\begin{itemize}
\item[] Corollary \ref{c:uniqueminframe} which says that $\mu ^0$ is
the unique (up to a homeomorphism) minimal measure function.
\item[] Corollary \ref{c:formofmeasure} which gives a partial
characterization of which continuous functions are realized as $\mu
(\fc)$ for some $\fc \in \fci$.
\end{itemize}

As has often been the case, the technique for proving these results is
to prove the corresponding result for sequences and sequence measure
functions and then apply the frame-sequence correspondence.

\subsection{Separable, irreducible and minimal sequence and  
frame measure functions} \label{subs:top1}

We begin by defining some natural classes of sequence and frame
measure functions.

\begin{Definition}
A sequence measure function $m:\cse\rightarrow\cstarw$ is 
\begin{itemize}
\item {\em separable} if
for every $v, \ w \in W$, $v \neq w$ there is $\xv\in\cse$ such that
$m(\xv)(v)\neq m(\xv)(w)$, 
\item{\em reducible} if there is a compact $V\subsetneq W$ such that
 $m':\cse\rightarrow\cc^*(V)$ is a sequence measure function, where
 $m'(\xv)=m(\xv){|}_V$,
\item {\em irreducible} if it is not reducible,
\item {\em minimal} if it is separable and 
irreducible.
\end{itemize}

\end{Definition}

\begin{Definition}
A frame measure function $m_f= m \circ b$ is ({\em separable,
reducible, irreducible, minimal}) if the corresponding sequence
measure function $m$ is (separable, reducible, irreducible, minimal).
\end{Definition}
The ultrafilter sequence and frame measure functions are not always
 {\em separable} as the following example shows:
\begin{Example}\label{exampl:6.3}
Suppose $I=\N$ and $I_n=\{1,2,\dots n \}$; therefore $|I_n|=n$.
Consider $p_1\in\N^*$ a free ultrafilter on $\N$, and define
\[ p_2 =\{ s+1~,~(s+1)\cup\{0\}~:~s\in p_1 \} \]
where $s+1=\{n+1: \ n\in s\}$

Notice $p_1\neq p_2$ since, for instance, $\{2k~:~k\in\N\}$ and
$\{2k+1~:~k\in\N\}$ would both be in $p_1$ which is impossible since
their intersection is empty.

For any $\xv \in X$ we have $x_n- x_{n-1} \leq |I_n| - |I_{n-1}| =1$
and $\frac{x_n}{|I_n|}= \frac{x_n}{n}\leq 1$.  Suppose $p_1- \lim
\frac{\xv}{\iv} =a$ thus for all $\epsilon >0$, there is a set $s \in
p_1$ for which $|\frac{x_n}{n} -a| < \epsilon$ for all $n\in s$.  Let
$N$ be such that $\frac{1}{N}< \epsilon$.  Note that $s'=\{ n \in s: n
\geq N\} \in p_1$ and set $t=s'+1 \in p_2$.  For $n\in t$,
\[ |\frac{x_n}{n} -a|=|\frac{x_n}{n} - \frac{x_{n-1}}{n-1} + \frac{x_{n-1}}{n-1} -a|= |\frac{x_n-x_{n-1}}{n} - \frac{x_{n-1}}{n(n-1)} + (\frac{x_{n-1}}{n-1} -a)|\]
\[ \leq \frac{1}{n} + \frac{1}{n} + \epsilon< 3 \epsilon . \]  
Thus $p_2 - \lim \frac{\xv}{\iv} =a$ as well, and so
$\mu(\cs)(p_2)=\mu(\cs)(p_1) $ . Therefore the set of continuous
functions $\mu(X)$ in $\cstarw$ does not separate $p_1$ from $p_2$ and
thus $\mu$ is not an example of a separable sequence measure function.
\end{Example}

We would like to use $\mu $ to construct a separable measure function.
Thus we are interested in grouping together all points in $\N ^*$
that produce the same values for all sequences.  To this end we
introduce the following equivalence relation on $\N^*$:
\begin{Definition}
For any $p_1,p_2\in\N^*$,
 we say $p_1\sim p_2$ if  $\mu(\xv)(p_1)=\mu(\xv)(p_2)$
 for all $\xv\in\cse$. 
\end{Definition}

It is easy to check that $\sim$ is an equivalence relation. Let $\No =
\N^*/\sim$. We consider $\No$ endowed with the {\em quotient
topology}: the finest topology such that the canonical projection
$\pi:\N^*\rightarrow\No$, $\pi(p)=\hat{p}=\{p'~|~p'\in\N^*~,~p'\sim
p\}$, is continuous. The open sets of $\No$ are therefore given by
$\{~U\subset\No : \pi^{-1}(U)~{\rm open~in }~\N^*\}$.

Considering $\No$ with the quotient topology we have:
\begin{itemize}

\item {\it $\No$ is compact} since it is the continuous image of the
compact space $\N^*$.

\item {The map \it $\mu ^0(\xv): \No \rightarrow \R$ defined by
$\mu^0(\xv)(\hat{p})=\mu(\xv)(p)$ is continuous for all $\xv \in
\cse$} since $\mu (\xv)$ is continuous on $\N^*$.
\item {\it $\No$ is Hausdorff} as the next two sentences show.  For
$p_1 \not= p_2 \in \No $, there must be a sequence $\xv$ for which
$\mu^0 (\xv)(p_1) \not= \mu^0 (\xv)(p_2)$ , and therefore there exist
disjoint open sets $U_1$, $U_2 \subset \R$ such that $\mu^0 (\xv)(p_1)
\in U_1$, $\mu^0 (\xv)(p_2)\in U_2$.  It follows that the open sets
$(\mu ^0 (\xv))^{-1}(U_1)$, $(\mu ^0 (\xv))^{-1}(U_2)$ separate $p_1$
and $p_2$.
\end{itemize}
The above allows us to define a new measure function:
\begin{Definition}
Denote by $\mu^0:\cse\rightarrow\cc^*(\No)$ the sequence measure
function defined as
\begin{equation}
\label{eq.m3.1a}
\mu^0(\xv)(\hat{p})=\mu(\xv)(p)
\end{equation}
Denote by $\mu ^0$ as well the corresponding frame measure function
$\mu^0 \circ b$.
\end{Definition}

We now show that $\mu^0 $ is minimal; in subsection \ref{subs:top2}
we will show that $\mu^0$ is essentially the unique minimal sequence
and frame measure function.  We begin by stating a trivial consequence
of Theorem \ref{t:canonembed}.
\begin{Corollary} \label{c:canonembed}
For any sequence measure function $m:\cse\rightarrow\cstarw$ there
exists an injection $\varphi : \No \rightarrow W$ such that
$m(\xv)(\varphi(p))= \mu^0 (\xv)(p)$ for all $\xv \in \cse$, $p \in
\No $.
\end{Corollary}
\pf The result follows trivially from Theorem \ref{t:canonembed} and
the definition of $\mu ^0$ which just eliminates the indistinguishable
points of $\N ^*$.

\begin{Proposition}
The map $\mu ^0: \cse \rightarrow\cc^*(\No)$ is a minimal sequence
measure function.
\end{Proposition}
\pf The definition of $\mu^0$ assures that it is separable.  Assume
that $\mu^0$ is not irreducible.  Thus there is a compact
$\N'\subsetneq \N^0$ so that $\mu ': \cse \rightarrow \cc^*(\N')$
defined by $\mu'(\xv)=\mu^0(\xv)|_{\N'}$ is again a sequence measure
function. Now consider a point $p \in \N^0 \backslash \N' $.  Denote
by $\varphi: \No \rightarrow \N'$ the map given in Corollary
\ref{c:canonembed}.  Thus $\mu'(\xv)(\varphi(p))= \mu(\xv)(p)$ for all
$\xv \in \cse$.  Since $\mu'(\xv)(\varphi(p))=\mu^0(\xv)(\varphi(p))$
for all $\xv \in \cse$, the separability of $\mu^0$ implies
$\varphi(p)=p$, a contradiction since $\varphi(p) \in \N'$, $p \in
\N^0 \backslash \N'$.  Thus $\mu^0$ must be irreducible and thus
minimal. \qed

As usual the above implies the corresponding result for frame measure
functions:

\begin{Corollary} 
The map $\mu^0: \fci \rightarrow\cc^*(\No)$ is a minimal frame
measure function.
\end{Corollary}

The construction above for getting $\No$ from $\N^*$ can be used for
any sequence or frame measure function $m:\cse\rightarrow\cstarw$ to
construct a separable sequence or frame measure function. Define on
$W$ the equivalence relation $v \sim w$ if $m(\xv)(v)=m(\xv)(w)$ for
all $\xv\in\cse$. The quotient space $W^0=W/\sim$ is then compact
Hausdorff with respect to the quotient topology.  We denote by $\pi$
the continuous map $\pi:W\rightarrow W^0$ defined by $\pi(v)=\pi(w)$
if and only if $v \sim w$. The sequence measure function $m$ induces a
map $m^0:\cse\rightarrow\cstarwo$ with
\begin{equation}
\label{eq.m3.1.1}
m^0(\xv)(p)=m(\xv)(q)~~,~~for~q\in\pi^{-1}(p).
\end{equation}
The definition of $m^0$ yields:
\begin{Proposition}\label{prop.m3.2.1}
The map $m^0:\cse\rightarrow\cstarwo$ is a separable sequence measure
function.  Consequently the map $m^0_f= m^0 \circ b$ that can be
constructed from a given frame measure function $m_f= m\circ b$ is a
separable frame measure function.
\end{Proposition}

\subsection{Uniqueness of the minimal sequence and frame measure function} 
\label{subs:top2}

\begin{Lemma}
\label{l.m3.51}
If $m:\cse\rightarrow\cstarw$ is minimal, then $\varphi :\No
\rightarrow W $ described in Corollary \ref{c:canonembed} is injective
with dense range.
\end{Lemma}

\pf Injectivity is a result of Corollary \ref{c:canonembed}.  If the range
$\varphi(\No)$ is not dense in $W$, then $m$ restricted to the closure
of $\varphi(\No)$ would also be a sequence measure function which
would contradict the minimality of $m$. \qed

\begin{Corollary}
\label{c.m2.c2}
For a minimal sequence measure function $m:\cse\rightarrow\cstarw$,
$m(\xv \wedge\yv)=min(m(\xv),m(\yv))$, $m(\xv \vee
\yv)=max(m(\xv),m(\yv))$ for any two sequences $\xv,\yv\in\cse$.
 
\end{Corollary}

\pf It follows from Proposition \ref{p:minmax} that the result is true for
the minimal sequence measure function $\mu^0$.  The result follows
from Lemma \ref{l.m3.51} and the continuity of the maps $m(\xv), \
m(\xv \wedge \yv)$ and $ m(\xv \vee \yv)$. \qed

\begin{Lemma}
\label{l.m.dens}
Let $m: \rs \rightarrow \cstarw$ be a minimal sequence measure
function.  For any $a,b\in \R$ and $v,w\in W$, there is an $\xv \in
\cse$ such that $m(\xv)(v)=a$ and $m(\xv)(w)=b$.
\end{Lemma}
\pf Recall $\iv =(|I_1|, |I_2|, \dots )$; $\iv$ is sequence compatible
and $m(\iv)(w)=1$ for all $w\in W$.

The case $a=b$ is simple since $m(a\iv)(w)=a$ for all $w\in W$. For
the case $a \not=b$, since $m$ is separable, there exists $\xv^0 \in
\cse$ such that $m(\xv^0)(v) \not= m(\xv^0)(w)$.  Let $c_1, \ c_2 \in
\R$ be determined by the linear system:
\[ c_1 m(\xv^0)(v) + c_2= a, \hspace{.5in} c_1 m(\xv^0)(w) + c_2=b.\]
Set $\xv= c_1 \xv^0 + c_2 \iv \in \rs$.  It follows by linearity of
the sequence measure function that $m(x)(v)=a$, $m(x)(w) =b$.  \qed

\begin{Theorem}[Density of Range]
\label{t.m.2.2}
Assume $m:\rs \rightarrow\cstarw$ is a minimal sequence measure function. Then
for every bounded real-valued continuous function $f\in\cstarw$, and every
$\eps>0$ there exists  $\xv\in\rs$  so that
$\norm{m(\xv)-f}_{\infty}<\eps$.
\end{Theorem}
\pf Lemma \ref{l.m.dens} coupled with the fact that $\cse$ is a
lattice with respect to $\vee$, $\wedge$ (Proposition \ref{p:lattice})
allows for the application of the lattice version of Stone's theorem
\cite{Naim72}, Chap. I, \S2,10.II ; the result is then immediate. \qed

\begin{Corollary}  
Given $m:\rs \rightarrow\cstarw$ a minimal sequence measure function,
for every real valued continous function $f\in \cstarw$ and every
$\eps >0$ there exists a constant $c$ and two frame compatible
sequences $\yv ^1$, $\yv ^2$, such that $\norm{c(m(\yv^1)-
m(\yv^2))-f}_{\infty} < \eps$.
\end{Corollary}

\pf Theorem \ref{t.m.2.2} establishes the existence of $\xv \in \rs$
for which $\norm{m(\xv)-f}_{\infty} < \eps$.  The result follows from
the fact that any $\xv \in \rs$ can be written as $\xv= c(\yv^1 - \yv
^2)$ with $\yv^1$, $\yv^2 $ frame compatible. \qed

As usual the above yields the corresponding result for frame measure
functions:

\begin{Corollary} \label{c:formofmeasure} 
Given $m:\fci \rightarrow\cstarw$ a minimal frame measure function,
for every real valued continous function $f\in \cstarw$ and every
$\eps >0$ there exists a constant $c$ and two frames $\fc ^1$, $\fc
^2$, such that $\norm{c(m(\fc^1)- m(\fc^2))-f}_{\infty} < \eps$.
\end{Corollary}

\begin{Lemma}
\label{l.m3.6}
If $m:\rs\rightarrow\cstarw$ is minimal, then $\varphi :\No
\rightarrow W $ described in Corollary \ref{c:canonembed} is
continuous.
\end{Lemma}

\pf To show continuity of $\varphi $ we will show that for all open
sets $V \subset W$ and all $p \in \varphi ^{-1}( V)$ there exists an
open set $U_p \in \No$ with $p \in U_p$ and $\varphi (U_p) \subset
V$. By Urysohn's Lemma, since $W\backslash V$ is closed, there is a
continuous function $\tilde{f}\in\cstarw$, so that $0\leq
\tilde{f}\leq 1$ on $W$, $\tilde{f}{|}_{W\setminus V}=1$, and
$\tilde{f}(\varphi (p))=0$.  By Theorem \ref{t.m.2.2} there exist
$\xv\in \rs$ such that $\norm{m(\xv)-\tilde{f}}_\infty\leq
\frac{1}{3}$.  Thus $m(\xv))|_{W \backslash V} \geq \frac{2}{3}$ and
$|m(\xv)(\varphi (p))| \leq \frac{1}{3}$.  Set $U_p= \mu^0 (\xv)
^{-1}( \ (- \frac{1}{2}, \frac{1}{2}) \ )$; $U_p$ is open (since $\mu
^0 (\xv)$ is continuous) and $p \in U_p$ (since
$|\mu^0(\xv)(p)|=|m(\xv)(\varphi (p))|=0\leq \frac{1}{3}$), and
$\varphi(U_p)= \mu^0(\xv)(U_p) \subset (- \frac{1}{2}, \frac{1}{2}) \
)$ whereas $m(x)(W\backslash V) \geq \frac{2}{3}$. \qed

\begin{Theorem} \label{t:uniquemin}
All minimal sequence measure functions $m:\rs\rightarrow\cstarw$are
topologically equivalent to $\mu ^0$, i.e. there exists a continuous
bijection with continuous inverse $\varphi: \No \rightarrow M$, such
that $m(\xv)(\varphi (p))= \mu ^0 (\xv)(p)$ for all $p \in \No$, $\xv
\in \rs$.
\end{Theorem}

\pf We let $\varphi: \No \rightarrow M $ be the map given in Corollary
\ref{c:canonembed}, Lemma \ref{l.m3.51}, and Lemma \ref{l.m3.6}.  From
these results we have that $\varphi $ is injective, has dense range,
and is continuous.  Since $\No$ is compact it follows from the
continuity of $\varphi$ that $\varphi (\No)$ is compact and thus it
must be all of $M$ (since it is dense in $M$).  Thus $\varphi $ is a
bijection.  Having established this bijection, we denote by $\varphi
^{-1}: M \rightarrow \No$ the inverse map.  The continuity of $\varphi
^{-1}$ is shown the same way as in Lemma \ref{l.m3.6}. \qed

\begin{Corollary} \label{c:uniqueminframe}
All minimal frame measure functions $m:\fci\rightarrow\cstarw$ are
topologically equivalent to $\mu ^0$, i.e. there exists a continuous,
bijection with continuous inverse $\varphi: \No \rightarrow M$, such
that $m(\fc)(\varphi (p))= \mu ^0 (\fc)(p)$ for all $p \in \No$, $\fc
\in \fci$.
\end{Corollary}

\begin{Remark}
We provide an example of a sequence measure function that is separable
but not minimal (that is it is not irreducible).  This implies the
existence of a frame measure function that is separable but not
minimal.  Let $|I_n|=2^n $ and consider the minimal measure function
$\mu^0 : \rs \rightarrow \cstarno $. Let $W= \No \cup \{w_0\}$ be 
the union of $\No$ with one extra point $w_0$.  
Pick two distinct $p_1$, $p_2 \in \N ^*$ so that $p_1$ contains the
set of odd integers, and $p_2$ contains the set of even integers.  
Define $m(\xv)(w_0)= \frac{1}{2}(\mu^0 (x) (p_1) + \mu^0(x) (p_2))$ and 
define $m(\xv) (p)= \mu^0 (x)(p)$.  Since $\No $ is a proper subset
of $W$, $m$ is not minimal. Now consider the frame compatible sequence $\xvt$
defined by
\[ \left\{ \begin{array}{rcl}
\tilde{x}_1 & = & 0 \\
\tilde{x}_{2n} & = & \tilde{x}_{2n-1} + |I_{2n}\setminus I_{2n-1}| \\
\tilde{x}_{2n+1} & = & \tilde{x}_{2n}
\end{array} \right. \]
Explicitly,
$\tilde{x}_{2n+1}=\tilde{x}_{2n}=\frac{2}{3}(4^n-1)$. Notice that
$\lim_{n\rightarrow\infty} \frac{\tilde{x}_{2n}}{|I_{2n}|} =
\frac{2}{3}$ whereas
$\lim_{n\rightarrow\infty}\frac{\tilde{x}_{2n+1}}{|I_{2n+1}|}
=\frac{1}{3}$. Now take a $p\in\No$. Then $m(\xvt)(p)$ equals either
$\frac{1}{3}$ or $\frac{2}{3}$ depending on whether $p$ contains the
set of odd integers, or not. In either case $\xvt$ separates $w_0$ from
$p$,
$$m(\xvt)(w_0)=\frac{1}{2}(\frac{1}{3}+\frac{2}{3}) = \frac{1}{2} \neq
m(\xvt)(p).$$
Thus $m$ is a separable but not minimal frame measure function.
\end{Remark}

\section{The \Cstar algebra of non-expansive operators}\label{sec4} 
\label{s:9}

Our approach to the classification of frames has been to examine the
sequence $b(\fc)$ associated to a frame $\fc$ via (\ref{eq:b}).  The
sequence $b(\fc)$ can be seen to be certain averages of the diagonal
elements of the Gram matrix $\{\ip{f_i}{\tilde{f}_j} \}_{i,j\in I}$.  We
now extend the definition of $b$ to all $I \times I$ matrices and then
compose this extended $b$ map with a sequence measure function $m$ to
give a measure on $I \times I$ matrices.  The result is an {\it
operator measure function} that resembles a trace on a large
subalgebra of operators.  In conjunction with some added structure on
the index set $I$, this expanded viewpoint leads to Theorem
\ref{t:redundancy} which states that  $m(\fc_1 \oplus
\fc_2)= m(\fc_1) + m(\fc_2)$ for a superframe $\fc_1 \oplus \fc_2$
where $\fc_1$ and $\fc_2$ need not be orthogonal but merely {\it
non-expansive} (see Definition \ref{d:framenonexp}).  This in turn
leads to a necessary density inequality for supersets of Gabor frames
(Theorem \ref{t9.3.2} and Corollary \ref{c9.4}).

We begin in Section \ref{s:7.1} by extending the definitions of
measure function and $b$ to the set of bounded operators.  We define
the important notion of {\it non-expansive} operators and frames and
show that the set of non-expansive operators is a large $C^*$
subalgebra of the set of bounded linear operators acting on $l^2(I)$.
We use this set up to prove the aforementioned result about supersets
in Section \ref{sec4.2}.

\subsection{Operator Measure Functions} \label{s:7.1}

We begin by defining $\cxs=\{x^1 + i x^2: x^1, \ x^2 \in \rs
\}$.  Recall the equivalence relation $\simm$ introduced in Definition
\ref{d:simm} applies to sequences in $\cxs$ as well. Thus $\xv \simm
\yv,$ $\xv, \yv \in \cxs $, if $\lim_{n\rightarrow \infty}
(x_n-y_n)/|I_n| =0$.

The following extends the map $b$ to operators.

\begin{Definition}
Let $b_{\mbox{op}}$ be the map from bounded linear operators on
$l^2(I)$ to sequences defined by
\[ 
b_{\mbox{op}} (A)= \{\sum_{i\in I_n} \ip{A \delta _i}{\delta_i} 
\}_{n \in \N}
\]
where $\{\delta_i\}_{i\in I}$ is the canonical basis of $l^2(I)$.
\end{Definition}

The range of $b_{\mbox{op}}$ lies in $\cxs$:

\begin{Proposition} \label{p:oper-sequence} For all $A \in B(l^2(I))$,
$b(A) \in \cxs$.
\end{Proposition}
\pf
Define 
\[a^+_j= \max( Re(\ip{A\delta _j}{\delta _j}),0), \hspace{.5in} \ a^-_j= 
\min( Re(\ip{A\delta _j}{\delta _j}),0),\]
\[ \ a^i_j= \max( Im(\ip{A\delta _j}{\delta _j}),0), \hspace{.5in} \  
a^{-i}_j= \min( Im(\ip{A\delta _j}{\delta _j}),0),\] thus $a^+_j +
a^-_j + i( a^i_j +a^{-i}_j) = \ip{A\delta _j}{\delta _j}$ with $a^+_j,
-a^-_j, a^i_j, -a^{-i}_j \leq ||A||$.  Define $x^+_n= \sum_{j \in I_n}
a^+_j$, $\ x^-_n= \sum_{j \in I_n} a^-_j$, $\ x^i_n= \sum_{j \in I_n}
a^i_j$, $\ x^{-i}_n= \sum_{j \in I_n} a^{-i}_j$.  It follows then that
$b_n(A) = x^+_n - (-x^- _n) + ix^i_n -i(-x^{-i})$.  It is
straightforward to verify that the sequences $\{x^+_n\}_{n \in \N}, \
\{-x^-_n\}_{n \in \N}, \ \{x^i_n\}_{n \in \N}, \ \{-x^{-i}_n\}_{n \in
\N} $ are all in $X^+$ (the appropriate $c$ being $||B||$) and thus
$b(A)= \{ b_n(A) \}_{n \in N} \in \cxs$. \qed

\begin{Remark} \label{r:framegram}
We note that given a frame $\fc$ and its associated Gram projection
$P\in B(l^2(I))$, we have $b(\fc)=b_{\mbox{op}}(P)$.
\end{Remark}

For the rest of this paper we will write $b$ for $b_{\mbox{op}}$.
Thus $b$ is both a map from frames to sequences (previous notation)
and the related map from linear operators to sequences.

Denote by $\cstarwcx$ the set of complex valued continuous maps on
$W$.  We now show that any sequence measure function has a unique
linear extension to $\cxs$.
\begin{Proposition} \label{p:cxextension}
Given a sequence function $m: \rs \rightarrow \cstarw $ , there exists
a unique linear map $\tm: \cxs \rightarrow \cstarwcx $ such that
$\tm|{\rs}=m$.
\end{Proposition}
\pf For any $\xv \in \cxs$, the decomposition of $\xv= \xv^1+ i \xv^2
$, $\xv^1$, $\xv^2 \in \rs$, is unique with $\xv^1_i = \mbox{Re}(\xv
_i)$, $\xv^2_i = \mbox{Im}(\xv_i)$.  Define $\tm= m(\xv^1) + i
m(\xv^2)$.  Thus $\tm$ is linear (since $m$ was linear) and $\tm
|_{\rs}$.  In addition $\tm$ is the unique linear extension since
there is only one way to write $\xv= \xv^1+ i \xv^2 $. \qed

 We now define an {\it operator measure } function :
 \begin{Definition}
 An operator measure function, $\bar{m}: B(l^2) \rightarrow \cxs $ is
 a map of the form $\bar{m}=\tm \circ b$ where $\tm$ is the linear
 extension of a sequence measure function described in Proposition
 \ref{p:cxextension}.
 \end{Definition}
 
We note that an operator measure function $\mbar $ is linear since it
is the composition of two linear maps.  The next few sections examine
the behaviour of $\mbar$.  We show that with added structure on the
index set $I$, there exists a large $C^*$ algebra $\cc \subset
B(l^2(I))$ for which $\mbar $ is tracial, i.e. $\mbar(AB)=\mbar(BA)$
for $A,B \in \cc$.  This tracial property is then used to prove
Theorem \ref{t3.3} which states that for a superframe $\fc_1\oplus\fc_2$ of
two non-expansive frames (see Definition \ref{d:framenonexp}) $\fc_1$,
$\fc_2$, the equation $m(\fc_1\oplus\fc_2) = m(\fc_1) + m(\fc_2)$
holds.


\subsection{The $C^*$ algebra of non-expansive operators\label{sec4.1}}

By a {\em quasi-distance} $d$ on $I$ we shall mean a map $d:I\times
I\rightarrow \R^+$ that satisfies: (i) $d(i,i)=0$, $d(i,j)\geq 0$;
(ii) $d(i,j)=d(j,i)$; (iii) $d(i,j)\leq d(i,k)+d(k,j)$, for any
$i,j,k\in I$.

For this section we shall consider an index set $I$ equipped with a
quasi-distance $d$. We call $(I,d)$ a {\em quasi-metric index set}.
We denote the ball of radius $R$ from $i\in I$ by
\begin{equation}
B_R(i)=\{j \in I: d(j,i) \leq R \}
\end{equation}

We shall say that $I$ has finite upper density with respect to $d$ if
$\sup _{i \in I}|B_R(i)| < \infty$ for all $R>0$.

Recall an algebra $S\subset \Bltwoi$ that is invariant under the
adjoint operation (i.e. $A^*\in S$ for any $A\in S$) is called a {\em
$C^{*}$ algebra} if it is closed in the operator norm topology.

\begin{Definition} \label{d:nonexpansive}
\begin{enumerate}
\item An operator $A \in B(l^2(I))$ is row {\em non-expansive} if for any $\epsilon >0$, there exists an $N(A,\eps) >0$ such that 
\begin{equation}
\label{eq.c.2.3}
\sum_{j\in I\backslash B_{N(A,\eps)}(i)}|\ip{A\delta_i}{\delta_j}|^2 < \eps
\end{equation}
for all $i\in I$.
\item An operator is non-expansive if both $A$ and $A^*$ are row non-expansive.  Denote by $\cc\subset B(l^2(I))$ the set of non-expansive operators.
\end{enumerate}
\end{Definition}

\begin{Theorem}
\label{t.c.1}
Suppose $I$ has finite upper density with respect to $d$.  Then $\cc$ is 
\begin{enumerate}
\item {\it closed under addition and scalar multiplication}, i.e. if
$A,B \in \cc$ and $c \in \C$ then $A+B \in \cc$ and $cA \in \cc$.
\item {\it closed under multiplication}, i.e. if $A,B \in \cc$ then
$AB \in \cc$.
\item{\it closed in the operator norm topology} i.e. given a filter
$\cal{J}$ on some set $S$ with $A_j \in \cc$ for all $j\in S$ and
$lim_{j\rightarrow \cal{J}} ||A-A_j||=0$ then $A\in \cc$.
\end{enumerate}
Consequently $\cc$ is a $C^{*}$ algebra.
\end{Theorem}

{\it Proof of 1.}

Fix an $\eps>0$. Set $N= \max(N(A, \frac{\eps}{4}), N(B,
\frac{\eps}{4}))$ with $N(A, \frac{\eps}{4})$, $N(B, \frac{\eps}{4})$
as in Definition \ref{d:nonexpansive}.  Thus for all $i$, we have
\[ \sum_{j \in I \backslash B_N(i)}|\ip{(A+B)\delta_i}{\delta_j}|^2 \leq 2(\sum_{j\in
I\backslash B_{N(A, \frac{\eps}{4})}(i)}|\ip{A\delta_i}{\delta_j}|^2 +
\sum_{j\in I\backslash
B_{N(B,\frac{\eps}{4})}(i)}|\ip{B\delta_i}{\delta_j}|^2) < \eps \]
This proves $A+B$ is non-expansive.

Setting $N=N(A,\frac{\eps}{|c|})$ yields 
\[ \sum_{j\in I\backslash B_N(i)}|\ip{cA\delta_i}{\delta_j}|^2 <\eps \]
for all $i\in I$, which proves $cA$ is non-expansive.

{\it Proof of 2.}

 Fix $\eps>0$.  Let $\eps_B=\frac{\eps}{4\norm{A}^2}$ and set
$N_B=N(B, \eps _B)$.  Let $\eps_A=\frac{\eps}{4 \norm{B}^2 D(N_B)}$,
where $D(N_B)=\sup_i |B_{N_B}(i)|$ (the upper bound on the number of
points of $I$ in a ball of radius $N_B$); set $N_A=N(A, \eps _A)$.
Let $N= N_A + N_B$ and fix $i \in I$.  We first note
\[ B\delta_i = \sum_{l\in I} \ip{B\delta_i}{\delta_l}\delta_l = v+ \sum_{l\in
 B_{N_B}(i)} 
\ip{B\delta_i}{\delta_l}\delta_l \]
for some vector $v$ with  $\norm{v}^2 <\eps_B$. Now
\begin{eqnarray*}
\sum_{j \in I \backslash B_N(i)}|\ip{AB\delta_i}{\delta_j}|^2 & = &
\sum_{j \in I \backslash B_N(i)}|\ip{Av}{\delta_j} + \sum_{l\in 
B_{N_B}(i)}\ip{B\delta_i}{\delta_l}\ip{A\delta_l}{\delta_j}|^2 
 \\ 
& \leq &  2\sum_{j\in I}
|\ip{Av}{\delta_j}|^2 + 2\sum_{j \in I \backslash B_N(i)}|\sum_{l\in B_{N_B}(i)}
\ip{B\delta_i}{\delta_l}\ip{A\delta_l}{\delta_j} |^2  \\ 
& \leq &  2\norm{A}^2 \eps_B + 2\sum_{j \in I \backslash B_N(i)} (D(N_B)\sum_{l\in B_{N_B}(i)}
|\ip{B\delta_i}{\delta_l}|^2 |\ip{A\delta_l}{\delta_j}|^2) 
\\ 
& = &   \frac{\eps}{2} + 2D(N_B)\sum_{l\in B_{N_B}(i)}|\ip{B\delta_i}{\delta_l}|^2
\sum_{j \in I \backslash B_N(i)} 
|\ip{A\delta_l}{\delta_j}|^2
\end{eqnarray*}
Now note that $I \backslash B_N(i) \subset I \backslash B_{N_A}(l)$ for any $l\in B_{N_B}(i)$.
Thus 
\[ \sum_{j \in I \backslash B_N(i)}|\ip{A\delta_l}{\delta_j}|^2\leq
\sum_{j \in I \backslash B_{N_A}(l)} 
|\ip{A\delta_l}{\delta_j}|^2 < \eps_A \]
and therefore:
\begin{eqnarray*}
\sum_{j \in I \backslash B_N(i)}|\ip{AB\delta_i}{\delta_j}|^2 & \leq &  \frac{\eps}{2} + 2D(N_B) \sum_{l\in
B_{N_B}(i)}|\ip{B\delta_i}{\delta_l}|^2 \eps_A  \\
& \leq& \frac{\eps}{2} + 2D(N_B) \norm{B\delta_i}^2 \eps_A  \\
& \leq & \frac{\eps}{2}+\frac{\eps}{2} = \eps
\end{eqnarray*}

{\it Proof of 3.}  Let $\eps>0$ be given. Then there is $K \in
\cal{J}$ so that for all $k\in K$,$A_k$ is non-expansive and
$\norm{A-A_k}^2 < \frac{\eps}{4}$. Let
$N_{\eps}=N(A_k,\frac{\eps}{4})$ for some fixed $k\in K$.  Then for
every $i \in I$,

\begin{eqnarray*}
\sum_{j \in I \backslash B_{N_{\eps}}(i)}|\ip{A\delta_i}{\delta_j}|^2 & = &
\sum_{j \in I \backslash B_{N_{\eps}}(i)}
|\ip{(A-A_k)\delta_i}{\delta_j} + \ip{A_k\delta_i}{\delta_j}|^2 \\
& \leq & 2 \sum_{j\in I} |\ip{(A-A_k)\delta_i}{\delta_j}|^2 + 2
\sum_{j \in I \backslash B_{N_{\eps}}(i)}
|\ip{A_k\delta_i}{\delta_j}|^2 \\
& \leq & 2\norm{A-A_j}^2 + \frac{\eps}{2} =\eps~~\qed
\end{eqnarray*}

\begin{Definition} \label{d:framenonexp}
We shall say that a frame $\fc$ is non-expansive if its associated
Gram projection is non-expansive.
\end{Definition}
Using elementary holomorphic functional calculus (see \S149 in
\cite{RieszNagy}) we can obtain the following:

\begin{Proposition}
\label{p.c.2.1} \label{p:polar}
Given a $C^*$ algebra $C$ acting on a Hilbert space and an operator
$A\in C$.  If the range of $A$ is closed then the orthogonal
projection onto the range of $A$ and the orthogonal projection onto
the range of $A^*$ are both in $C$.
\end{Proposition}

This result has a couple of consequences: it gives a simpler
sufficient (but not necessary) condition for a frame to be
non-expansive (Corollary \ref{c.c.1} below) and it plays a key role in
the proof of Theorem \ref{t3.3}.

\begin{Corollary}\label{c.c.1}
For any frame $\fc\in\fci$, if its Gram operator
$G:\ltwoi\rightarrow\ltwoi$, $G(c) =\{\sum_{j\in
I}\ip{f_j}{f_i}c_j{\}}_{i\in I}$ is non-expansive, then the $\fc$ is
non-expansive, as are the associated Parseval frame and the canonical
dual frame.
\end{Corollary}

\pf If $G$ is non-expansive, $G\in\cc$. Since $\fc$ is frame, the range of
$G$ is closed. Thus the associated Gram projection, by Proposition
\ref{p:polar}, is also in $\cc$, and thus $\fc$ is non-expansive.
Since $\fc$, the associated Parseval frame
$\fc^{\#}=\{S^{-1/2}f_i\}$ and the canonical dual frame
$\tilde{\fc}=\{S^{-1}f_i\}$ all have the same associated Gram
projection, they are all non-expansive. \qed

\begin{Remark}

Corollary \ref{c.c.1} is merely a sufficient condition as the
following construction demonstrates.  Let $S$ be a self-adjoint
operator that is not non-expansive.  It follows that the invertible
operator $G=S +2\norm{S}I$ is also not non-expansive. In this case,
the frame $\gc=\{g_i=G^{1/2}\delta_i\}$ is a Riesz basis and hence is
non-expansive (since the corresponding projection for a Riesz basis is
the identity).  However, the frame $\gc$ has a non-expansive Gram
operator $G$.
 
\end{Remark}

\subsection{The measure function and supersets\label{sec4.2}}

In this subsection we show that condition 4. of Definition
\ref{d:sequencemf} can be extended to non-orthogonal superframes that
are non-expansive. In particular we obtain a density-type result.

The main result that allows us to develop the theory is the tracial
property of the extended measure $\mbar$ on $\cc$ (Lemma \ref{l.c.4}).
The result will hold when the quasi distance $d$ and the decomposition
$I=\cup_nI_n$ have the following compatibility which essentially says
that the boundary (with respect to $d$) of subsets $(I_n)_{n\geq 0}$
are asymptotically smaller than their interior:

\begin{Definition}
The collection $(I,d,(I_n)_n)$ is called a  
{\em \metricindexset} if the quasi distance $d$ has finite upper density and 
for all $R>0$,
\begin{equation}
\label{eq.c.2.2bis}
\ \lim_{n\rightarrow\infty} \frac{|\cup_{j\in I\backslash I_n} B_R(j)\cap
I_n|}{|I_n|} = 0
\end{equation}\end{Definition}

\begin{Lemma}\label{l.c.4} Assume $(I,d,(I_n)_n)$ is a~ \metricindexset.
Then for any two non-expansive operators $T_1,T_2\in\cc$,
\begin{equation}
\label{eq.c.2.1}
\mbar(T_1T_2)=\mbar(T_2T_1)
\end{equation}
\end{Lemma}

\pf

Equation \ref{eq.c.2.1} is equivalent to 
\begin{equation}
\label{eq.c.2.2}
\lim_{n\rightarrow \infty} \frac{1}{|I_n|} b_n(T_1T_2 - T_2T_1)=
\lim_{n\rightarrow \infty} \frac{1}{|I_n|} (b_n(T_1T_2)- b_n(T_2T_1))=0
\end{equation}
Recall that $T\in\cc$ implies that both $T$ and $T^*$ are non-expansive.  
Since  $\{\delta_i\}_{i\in I}$ is an orthonormal basis:
\[ \frac{1}{|I_n|} b_n(T_1T_2) = \frac{1}{|I_n|}\sum_{i\in I_n} \sum_{j\in I}
\ip{T_2\delta_i}{\delta_j}\ip{T_1\delta_j}{\delta_i} \] 
Using the
corresponding expansion for $\frac{1}{|I_n|} b_n(T_2T_1)$ and
subtracting from the above, we get
\begin{equation}
\label{eq.c.2.3a}
\frac{1}{|I_n|}b_n(T_1T_2-T_2T_1) = \frac{1}{|I_n|}\sum_{i\in I_n}
\sum_{j\not\in I_n}
\ip{T_1\delta_j}{\delta_i}\ip{T_2\delta_i}{\delta_j} - \frac{1}{|I_n|}
\sum_{i\not\in I_n}\sum_{j\in I_n}
\ip{T_1\delta_j}{\delta_i}\ip{T_2\delta_i}{\delta_j}
\end{equation}
We shall show that the right hand side of (\ref{eq.c.2.3a}) has limit
$0$ as $n\rightarrow \infty$ which will establish the result.  We
apply Cauchy-Schwarz to the first term on the right side of
(\ref{eq.c.2.3a}) and obtain
\begin{eqnarray}\label{eq.c.2.3aa}
  &  & |\frac{1}{|I_n|}\sum_{i\in I_n}\sum_{j\not\in I_n}
\ip{T_1\delta_j}{\delta_i}\ip{T_2\delta_i}{\delta_j}|^2 \\
 & \leq &
(\frac{1}{|I_n|}\sum_{i\in I_n}\sum_{j\in I\setminus I_n} |\ip{T_1^*\delta_i}{
\delta_j}|^2)
(\frac{1}{|I_n|}\sum_{i\in I_n}\sum_{j\in I\setminus I_n} |\ip{T_2\delta_i}{
\delta_j}|^2)
\end{eqnarray}
Fix $\eps>0$. Let $N$ be a radius in the definition of
non-expansiveness that works for $T_1,T_2,T_1^*,T_2^*$ simultaneously.
Write $I_n=J_n\cup D_n$ where $D_n=I_n\cap(\cup_{j\in I\setminus I_n}
B_N(j))$ is the set of points of $I_n$ that are within distance $N$ of
the boundary, and $J_n=I_n\setminus D_n$ is the rest. Decomposing the
sums over $i\in I_n$ into the sums over $D_n$ and $J_n$, we have that
(\ref{eq.c.2.3aa}) is bounded above by

\begin{eqnarray}
 &   & (\eps + \frac{1}{|I_n|}\sum_{i\in D_n} 
\sum_{j\in I}||\ip{T_1^*\delta_i}{\delta_j}|^2)(\eps +
\frac{1}{|I_n|}\sum_{i\in D_n}  \sum_{j\in I}||\ip{T_2\delta_i}{\delta_j}|^2)
\nonumber \\
 & \leq & (\eps + \frac{|D_n|}{|I_n|} \norm{T_1}^2)(\eps + \frac{|D_n|}{|I_n|}
\norm{T_2}^2) \nonumber
\end{eqnarray}
A similar inequality is obtained for the second term in
(\ref{eq.c.2.3a}) and thus
\[ |\frac{1}{|I_n|}b_n(T_1T_2-T_2T_1)| \leq 2(\eps+\frac{|D_n|}{|I_n|}A) \]
where $A=max(\norm{T_1}^2,\norm{T_2}^2)$. Using the asymptotic assumption
(\ref{eq.c.2.2bis}) we obtain
\[ \lim_{n\rightarrow\infty} |\frac{1}{|I_n|}b_n(T_1T_2-T_2T_1)| \leq 3\eps \]
Since $\eps$ was arbitrary, we obtain (\ref{eq.c.2.2}). \qed

We now prove that frame measure functions are linear on supersets of
non-expansive frames:

\begin{Theorem}
\label{t3.3} \label{t:redundancy}
Assume $(I,d,(I_n)_n)$ is a~ \metricindexset~ and
 $m:\fci\rightarrow\cstarm$ a frame measure function. Suppose
 $(\fc_1,\fc_2)$ is a superframe of two non-expansive frames. Then
 $\fc_1\oplus\fc_2$ is non-expansive and
\begin{equation}
\label{eq.c.2.7}
m(\fc_1\oplus\fc_2) = m(\fc_1) + m(\fc_2)
\end{equation}
\end{Theorem}

\pf
We first show that $\fc_1\oplus\fc_2$ is non-expansive. Let $P_1,P_2$
denote the associated Gram projections to the two frames $\fc_1$ and
$\fc_2$. The definition of non-expansive frames gives $P_1,P_2\in\cc$.
Since $\fc_1\oplus\fc_2$ is a frame, we have by Proposition \ref{p:A2}
that $P_1 +P_2$ has closed range and thus by Proposition
\ref{p:polar}, the projection onto the range of $P_1 + P_2$, which is
the associated Gram projection for $\fc_1 \oplus \fc_2$, is also
non-expansive.

Let $P$ be the associated Gram projection for $\fc_1\oplus\fc_2$,
i.e. $P$ is the projection onto the range of $P_1+P_2$, that is
$P=P_1\vee P_2$. The statement (\ref{eq.c.2.7}) is equivalent to
proving
\begin{equation}
\label{eq.c.2.9}
 \mt(P) = \mt(P_1)+\mt(P_2)
\end{equation}

Consider $A=P_1-P_1P_2$. The superframe condition amounts (equivalently) to the
condition that $\norm{P_1P_2}<1$. Hence, when restricted to $Ran\,P_1$, 
$A=1-P_1P_2$ is invertible, hence its range is $Ran\,P_1$. Therefore $Ran\,A$
is closed, and equals $Ran\,P_1$. On the other hand any $x\in\ltwoi$ admits a
unique decomposition $x=x_1+x_2+x'$, where $x_1\in Ran\,P_1$, $x_2\in
Ran\,P_2$, and $x'\in Ran\,(1-P)$. Then
$\norm{Ax}=\norm{Ax_1} \geq (1-\norm{P_1P_2})\norm{x_1}$. Hence 
$ker\,A=ker\,(P-P_2)$ which implies $(ker\,A)^{\perp}=Ran\,(P-P_2)$. 
Since $A$ is in $\cc$ the partial isometry $V$ of the polar
decomposition $A=V(A^*A)^{1/2}$ belongs to $\cc$ using again standard 
holomorphic functional calculus arguments (as in \cite{RieszNagy}). 
Furthermore $V$ has initial space $Ran\,(P-P_2)$,
and final space $Ran\,P_1$, that is $VV^* = P_1$, and $V^*V=P-P_2$.
Since $\mt$ is tracial on $\cc$, it follows
$\mt(P_1)=\mt(VV^*)=\mt(V^*V)=\mt(P-P_2)$. But $P=(P-P_2)+P_2$ is an
orthogonal decomposition of $P$, therefore $\mt(P)=\mt(P-P_2)+\mt(P_2)$, which
together with the previous relation proves (\ref{eq.c.2.9}) and the
Theorem. \qed

The following corollary  immediately follows using induction:

\begin{Corollary}\label{c3.3}
Assume $(\fc_1,\cdots,\fc_D)$ is a superframe of non-expansive frames. Then
$\fc_1\oplus\cdots\oplus\fc_D$ is non-expansive and
\begin{equation}
\label{eq.c.2.8}
m(\fc_1\oplus\cdots\oplus\fc_D)=m(\fc_1)+\cdots+m(\fc_D)
\end{equation}
\end{Corollary}

\section{Measure functions and the index set} \label{s:index}

In this section we study how different frame indexing and finite
averaging methods affect the measure function and the property of
non-expansiveness. Because all measure functions contain a copy of the
ultrafilter measure function $\mu$ (cf Corollary \ref{c:frameembed})
we shall consider only the case of the ultrafilter frame measure
function $\mu$, and comment on the extension of these results to
arbitrary frame measure functions.

Assume $I$ and $J$ are countable index sets, and $a:I\rightarrow J$ is
a bijection. Assume also $(I_n)_n$ and $(J_n)_n$ are nested sequences
of finite subsets covering $I$, respectively $J$.  Our goal is to
establish how equivalence classes of frames in $\fc[I]$ are related to
equivalence classes of frames in $\fc[J]$. More generally, we will
examine the correspondence of operators between $B(l^2(I))$ and
$B(l^2(J))$ and the preservation of the non-expansiveness property.

First we note that the map $a$ induces a mapping on frames:
\begin{equation}
a_{*}:\fc[J]\rightarrow\fc[I] \label{eq:8.1.1}
~~,~~a_{*}(\fc)=\{f_{a(i)}~;~i\in I\}
\end{equation}
and a mapping on operators:
\begin{equation}
a_{*}:B(l^2(J))\rightarrow B(l^2(I))~~,~~
\ip{a_{*}(T)\delta_{i_1}}{\delta_{i_2}}
= \ip{T\epsilon_{a(i_1)}}{\epsilon_{a(i_2)}} 
\label{eq:8.1.2}
\end{equation}
where $(\delta_i)_i$ and $(\epsilon_j)_j$ are the canonical bases of
$l^2(I)$ and $l^2(J)$ respectively. 

We are interested in the following tasks:

\begin{enumerate}
\item {\it Measure Preservation.}  Find conditions on $a$ so that for
all operators $T\in B(l^2(J))$, the ultrafilter frame measure
functions for $T$ and $a_{*}(T)$ are equal.

\item {\it Non-expansiveness Preservation.}

Assuming that $(I,d)$ and $(J,e)$ are quasi-metric index sets, find
conditions on $a$ so that for all operators $T\in B(l^2(J))$, $T$ is
non-expansive if and only if $a_{*}(T)$ is non-expansive.  In
particular we obtain that $\fc\in\fc[J]$ is non-expansive if and only
if $a_{*}(\fc)$ is non-expansive.
\end{enumerate} 

We address each of these in the subsequent two sections.

\subsection{Measure preserving indexing}

The following gives a condition for $a$ that preserves the value of
the measure function.

\begin{Proposition}\label{prop8.1}
If the map $a:I\rightarrow J$ satisfies the following property
\begin{equation}
\label{eq:8.2.1}
\lim_n\frac{|a(I_n)\cap J_n|}{|I_n|}=\lim_n \frac{|J_n|}{|I_n|}=1
\end{equation}
then $\mu(T)=\mu(a_{*}(T))$ for all $T\in B(l^2(J))$. Explicitely
this means:
\begin{equation}
\label{eq:8.2.3}
\plim_{n}\frac{1}{|J_n|}\sum_{j\in J_n}T_{j,j}
=\plim_{n}\frac{1}{|I_n|}\sum_{i\in I_n}T_{a(i),a(i)}
\end{equation}
for all $p\in\N^*$.
\end{Proposition}

\pf

Since $T$ is bounded, it follows $|T_{j,j}|\leq r:=\norm{T}$ for all $j$.
First we have:
\begin{equation}
\frac{1}{|J_n|}\sum_{j\in J_n}T_{j,j}-\frac{1}{|I_n|}\sum_{i\in I_n}
T_{a(i),a(i)} = \frac{1}{|J_n|}\sum_{j\in J_n\setminus a(I_n)} T_{j,j}
 + \frac{|I_n|-|J_n|}{|I_n|\cdot|J_n|}\sum_{j\in J_n\cap a(I_n)}T_{j,j}
-\frac{1}{|I_n|}\sum_{j\in a(I_n)\setminus J_n}T_{j,j}
\end{equation}
Upper bounding each term, we get:
\begin{equation}
|\frac{1}{|J_n|}\sum_{j\in J_n}T_{j,j}-\frac{1}{|I_n|}\sum_{i\in I_n}
T_{a(i),a(i)}|\leq r\frac{|J_n\setminus a(I_n)|}{|J_n|}+r\frac{||I_n|-|J_n||
\cdot|J_n\cap a(I_n)|}{|I_n|\cdot|J_n|}+r\frac{|a(I_n)\setminus J_n|}{|I_n|}
\end{equation}
Condition (\ref{eq:8.2.1}) implies now that each term tends to zero as $n$
goes to infinity. Hence we get:
\[ lim_n [\frac{1}{|J_n|}\sum_{j\in J_n}T_{j,j}-\frac{1}{|I_n|}\sum_{i\in I_n}
T_{a(i),a(i)}] = 0 \]
which implies (\ref{eq:8.2.3}). \qed

\begin{Remark}
The same condition (\ref{eq:8.2.3}) guarantees the preservation
of equivalence classes of frames, that is
for all $\fc^1,\fc^2\in\fc[J]$
$\fc^1\simm_{J}\fc^2$ if and only if $a_{*}(\fc^1)\simm_{I} a_{*}(\fc^2)$.

Thus, in general, an arbitrary frame measure function on $\fc[I]$, 
$m:\fc[I]\rightarrow \cstarm$,  induces a measure function on $\fc[J]$, 
$a^{*}(m):\fc[J]\rightarrow \cstarm$ via $a^{*}(m)(\fc)=m(a_{*}(\fc))$.
\end{Remark}

\subsection{Indexing preserving non-expansiveness}

Now we examine when non-expansive operators are pulledback through
$a_{*}$ into non-expansive operators.  We use the same setting as
before where now $(I,d)$ and $(J,e)$ are assumed to be quasi-metric
index sets and $a:I\rightarrow J$ is the bijection.  We have the
following result:
\begin{Proposition}\label{prop8.2}
Suppose there exists a function $r:[0,\infty) \rightarrow [0, \infty)$ 
such that 
\begin{equation}
\label{eq:8.3.1}
\forall j_1,j_2\in J~~~d(a^{-1}(j_1), a^{-1}(j_2)) < r(e(j_1, j_2))
\end{equation}
Then if $T\in B(l^2(J))$ is non-expansive, then $a_{*}(T)$ is non-expansive in
 $B(l^2(I))$.
\end{Proposition}

\pf

Assume  that $T$ is non-expansive and choose an arbitrary $\eps>0$. 
Set $N=N_\eps$ from the non-expansive definition for $T$, then:
\begin{eqnarray}
 \sum_{i'\in I, d(i,i')>r(N)}|\ip{a_{*}(T)\delta_i}{\delta_{i'}}|^2
 & = & \sum_{j'\in J,d(i,a^{-1}(j))>r(N)}|\ip{T\epsilon_{a(i)}}{
\epsilon_{j'}}|^2 \nonumber \\
& \leq & \sum_{j'\in J,e(a(i),j')>N}|\ip{T\epsilon_{a(i)}}{\epsilon_{j'}}|^2
<\eps. \nonumber
\end{eqnarray}
A similar argument holds  for $T^*$ and thus $a_{*}(T)$ is non-expansive. \qed

\begin{Remark}
An immediate consequence of this result is that if $\fc\in\fc[J]$
is non-expansive then $a_{*}(\fc)$ is non-expansive as well.
\end{Remark}

\begin{Remark}
If the two quasi-metric spaces $(I,d)$ and $(J,e)$ satisfy the
assumption of Proposition \ref{prop8.2}, then one can always choose a
continuous and monotonically inreasing $r$ in (\ref{eq:8.3.1}).
\end{Remark}

\subsection{A Consequence}

Now we can put together Theorem \ref{t3.3}, and Propositions 
\ref{prop8.1}, \ref{prop8.2}, and obtain the following

\begin{Theorem}
\label{t8.1}
Assume $(I,d,(I_n)_n)$ is a \metricindexset~ and $(J,e,(J_n)_n)$
is so that $(J,e)$ is a quasi-metric index set. Assume $a:I\rightarrow J$
is a bijection that satisfies
\begin{equation}
\label{eq:8.4.1}
\lim_n\frac{|a(I_n)\cap J_{n}|}{|I_n|}=\lim_n 
\frac{|J_{n}|}{|I_n|}=1 
\end{equation}
and there exists a function $r:[0,\infty) \rightarrow [0, \infty)$ such that 
\begin{equation}
\label{eq:8.4.2}
\forall j_1,j_2\in J~~~d(a^{-1}(j_1), a^{-1}(j_2)) < r(e(j_1, j_2))
\end{equation}

Assume $\fc^1\in\fc[I]$ is non-expansive with respect to the quasi-metric 
index set $(I,d)$ and $\fc^2\in\fc[J]$ is non-expansive with respect to the
quasi-metric index set $(J,e)$. Then, if 
$\fc=\{f^1_i\oplus f^2_{a(i)}~;~i\in I\}$
is frame (that is, $(\fc^1,a_{*}(\fc^2))$ is a superframe) then
$\fc$ is nonexpansive with respect to $(I,d)$ and
\begin{equation}
\label{eq:8.4.3}
\mu(\fc)(p)=\mu(\fc^1)(p)+\mu(\fc^2)(p)~~,~~\forall p\in\N^*.
\end{equation}
Explicitly, for every free ultrafilter $p\in\N^*$,
\begin{equation}
\label{eq:8.4.4}
\mu(\fc)(p) = \plim \frac{1}{|I_n|}\sum_{i\in I_n}\ip{f^1_i}{\tilde{f^1_i}}
+\plim \frac{1}{|J_n|}\sum_{j\in J_n}\ip{f^2_j}{\tilde{f^2_j}}.
\end{equation}
\end{Theorem}

This statement can be straightforwardly extended to a finite
collection of frames that form a superframe.

One can replace the free ultrafilter frame measure function $\mu$ by
any other frame measure function $m$ on $\fc[I]$; consequently, in
this case we have: 
\begin{equation}
\label{eq:8.4.5}
m(\fc)(x) = m(\fc^1)(x) + a^{*}(m)(\fc^2)(x).
\end{equation}

\section{Application to Gabor Frames and Superframes} \label{s:8} \label{s:gab}

In this section, we apply our results to Gabor frames and superframes.
We begin with some added notation and preliminaries.

For a function $g\in L^2(\R^m)$, a point $\lambda=(t,\omega)\in
\R^m \times \R^m$, and a phase $\varphi_\lambda\in\R$
 denote by $g_\lambda(x)=e^{i\varphi_\lambda}e^{2\pi
i\ip{\omega}{x}}g(x-t)$ the $\lambda$-time-frequency shift of $g$. 

\begin{Definition}
Given a function $g \in L^2(\R^m)$ and a set of time-frequency shifts
$\Lambda \subset \R^m \times \R^m$, and phases
$\{\varphi_\lambda\}_{\lambda\in\Lambda}$ define the {\it Gabor set}
$\G (g, \Lambda)=\{g_{\lambda} \}_{\lambda \in \Lambda}$.  A {\it
Gabor frame} is a Gabor set that is a frame sequence.
\end{Definition}

For ease of notation we will omit the explicit mention of the phase
system $\{\varphi_\lambda\}_\lambda$. 

We define $Q_n(c)=\{\lambda\in\R^{2m}~|~\norm{\lambda-c}_\infty\leq
\frac{n}{2}\}$ to be the box inside $\R^{m}\times\R^m$ centered at
$c\in\R^{2m}$ and of size length $n$. 

Given a Gabor set $\G(g,\Lambda)$, the most natural way of indexing is
given by the set $\Lambda$ itself. Thus
 $(\Lambda,\norm{\cdot}_\infty)$ becomes a quasi-metric index set. 
Note that $\norm{\cdot}_\infty$ may not be a
distance because we allow repetitions of the same time-frequency point
in $\Lambda$. 

We need to define the nested sequence of finite subsets
$(\Lambda_n)_n$.  Fix a center $O\in\R^{2m}$ (not necessarily the
origin).  It turns out that the natural choice of
$\Lambda_n=Q_n(O)\cap\Lambda$ is not suitable for measuring Gabor
frames. To fix this issue we instead replace $Q_n(O)$ by a ``skewed''
tile $MQ_n(O)$, where $M$ is a suitable $2m\times 2m$ invertible
matrix. We can do this either by simply defining
$\Lambda_n=(MQ_n(O))\cap\Lambda$, or by changing the distance in
$\R^{2m}$ and replacing $\norm{x}_\infty$ by
$\norm{x}_{M,\infty}:=\norm{M^{-1}x}_{\infty}$. The two approaches are
equivalent. However for simplicity of computations we will adopt the
former approach, namely we keep the $\norm{}_\infty$ distance in
$\R^{2m}$ and define $\Lambda_n=(MQ_n(O))\cap\Lambda$.

We will compute the free ultrafilter frame measure function of
$\G(g,\Lambda)$ with respect to partition $(\Lambda_n)_n$.  We will
show that $(\Lambda,\norm{\cdot}_\infty,(\Lambda_n)_n)$ is a
~\metricindexset~, and $\G(g,\Lambda)$ is non-expansive. Next we compute
the frame measure function from Gabor superframes and obtain a necessary
density type condition.

\subsection{Free ultrafilter frame measure function of Gabor frames}
\label{s:8.1}

Let us consider a Gabor frame $\G(g,\Lambda)$. Then the upper and
lower Beurling densities of $\Lambda$, $D^{+}_B(\Lambda)$, and
$D^{-}_B(\Lambda)$, satisfy (see the historical note \cite{he06-1} of
this result)
\[ 1\leq D^{-}_B(\Lambda) \leq D^{+}_B(\Lambda) <\infty \]
where 
$$D^{+}_B(\Lambda) = \limsup_{n}\sup_{c\in\R^{2m}} \frac{|\Lambda\cap
Q_n(c)|}{n^{2m}}~,~D^{-}_B(\Lambda)=\liminf_n
\inf_{c\in\R^{2m}}\frac{|\Lambda\cap Q_n(c)|}{n^{2m}}.$$ In particular
this means there is a size $L_0>0$ and an integer $U_0\geq 1$ so that
every box of side length $L_0$ in $\R^{2m}$ contains at least one
point of $\Lambda$ and at most $U_0$ points of $\Lambda$. Fix a point
$O\in\R^{2m}$, an invertible matrix $M$ in $\R^{2m\times 2m}$ and let
$\Lambda_n=\Lambda\cap MQ_n(O)$ as before. For any length $R$, the box
$Q_n(O)$ is covered by at most $(\frac{n}{R}+1)^{2m}$ boxes of side
length $R$, and includes at least $(\frac{n}{R}-1)^{2m}$ disjoint
boxes of side length $R$.  For the skewed box $MQ_n(O)$ the situation
is the following. There are two numbers $c_1(M)$ and $c_2(M)$
depending on the matrix $M$ so that, at most
$c_1(M)(\frac{n}{R})^{2m}+c_2(M)(\frac{n}{R})^{2m-1}$ boxes are needed
to cover $MQ_n(O)$, and at least
$c_1(M)(\frac{n}{R})^{2m}-c_2(M)(\frac{n}{R})^{2m-1}$ disjoint boxes
of side length $R$ are included inside $MQ_n(O)$.
 With this set up we have the following:
\begin{Theorem}
\label{t:umis}
The collection $(\Lambda,\norm{\cdot}_\infty,(\Lambda_n)_n)$ is a~ 
\metricindexset.
\end{Theorem}

\pf $(\Lambda,\norm{\cdot}_\infty)$ has finite upper density since
every ball of radius $R$ contains at most $(\frac{2R}{L_0}+1)^{2m}$
boxes of side length $L_0$, and every box of side length $L_0$ has at
most $U_0$ points.  The second condition (\ref{eq.c.2.2bis}) 
is proved as follows. On the one hand for large $n$, each $\Lambda_n$ has
the cardinal bounded by:
\[ c_1(M)(\frac{n}{L_0})^{2m}-c_2(M)(\frac{n}{L_0})^{2m-1} 
\leq |\Lambda_n| \leq \left( c_1(M)(\frac{n}{L_0})^{2m}+c_2(M)
(\frac{n}{L_0})^{2m-1} \right)U_0 \]
On the other hand 
$$\cup_{j\in\Lambda\setminus \Lambda_n} B_R(j)\cap\Lambda_n
= \left( M(Q_n(O)\setminus Q_{n-R}(O)) \right) \cap \Lambda $$
Hence
\begin{eqnarray}\nonumber
&  |\cup_{j\in\Lambda\setminus\Lambda_n} B_R(j)\cap\Lambda_n | \leq & \\
& ((c_1(M)(\frac{n}{L_0})^{2m}+c_2(M)(\frac{n}{L_0})^{2m-1})
 -(c_1(M)
(\frac{n-2R}{L_0})^{2m}-c_2(M)(\frac{n-2R}{L_0})^{2m-1}))U_0 & \nonumber
\end{eqnarray}
Putting these two estimates together we obtain
\[ \lim_{n\rightarrow\infty} \frac{|\cup_{j\in\Lambda\setminus\Lambda_n}
B_R(j)\cap\Lambda_n|}{|\Lambda_n|} = 0. \qed \]

Consider a Gabor frame $\G(g,\Lambda)$ for $L^2(\R^m)$. 
Fix a point $O\in\R^{2m}$, an invertible matrix $M\in\R^{2m\times 2m}$,
 and set $\Lambda_n=\Lambda\cap MQ_n(O)$ as before.
For any free ultrafilter $p\in\N^*$, the set $\Lambda$ has density:
\begin{equation}
\label{eq:dens}
D(\Lambda;p,M) = \plim \frac{|\Lambda_n|}{vol(MQ_n(O))}
=\plim\frac{|\Lambda\cap(MQ_n(O))|}{det(M) n^{2m}}
\end{equation}

We recall a fundamental result obtained in \cite{bacahela06,bacahela06-1}.

\begin{Theorem}[\cite{bacahela06}]
\label{t:bchl}
Assume $\G(g,\Lambda)$ is a frame for $L^2(\R^m)$ and $\{\tilde{g}_\lambda~;~
\lambda\in\Lambda\}$ is its canonical dual frame. Then for any free ultrafilter
$p\in\N^*$,
\begin{equation}
\label{eq:bchl}
\plim \frac{1}{|\Lambda_n|}\sum_{\lambda\in\Lambda_n}
\ip{g_\lambda}{\tilde{g}_\lambda} = \frac{1}{D(\Lambda;p,M)}
\end{equation}
\end{Theorem}

The fact that we use skewed boxes instead in regular boxes does not affect 
the result. As we mentioned earlier, we can change the metric to account
for the skewness, and apply directly the results from \cite{bacahela06,bacahela06-1}.

This fundamental relation gives us a simple way to compute the
free ultrafilter frame measure function of irregular Gabor frames
 (compare to Theorem 3 in \cite{bacahela06-1}):

\begin{Theorem}\label{t:gabormeasure}
For any Gabor frame $\G(g,\Lambda)$ and indexing $(\Lambda,(\Lambda_n)_n)$
as before, the free ultrafilter frame measure function is
\begin{equation}
\label{eq:gm}
\mu(\G)(p)=\frac{1}{D(\Lambda;p,M)}~~,~~\forall p\in\N^*
\end{equation}
\end{Theorem}

\begin{Remark}
If $\Lambda$ has uniform density $D_0$ (that is
$D^{-}_B(\Lambda)=D^{+}(\Lambda)=D_0$) then
$\mu(\G)=\frac{1}{D_0}1_{\N^*}$, that is, the measure function of $\G$
is the constant function $\frac{1}{D_0}$, independent of the matrix
$M$.  In fact, for any measure function $m:\fc[\Lambda]\rightarrow
\cstarw$ the measure of $\G$ is $m(\G)=\frac{1}{D_0}1_W$.

For $\Lambda=A\Z^{2m}$ for some invertible matrix $A$, then
 $D_0=\frac{1}{det(A)}$ regardless of matrix $M$,
 and thus $m(\G)=(det(A))1_W$. In particular,
 for $\Lambda=\al\Z^m\times\be\Z^m$, $D_0=\frac{1}{(\al\be)^{m}}$
 and $m(\G)=(\al\be)^m 1_W$.
\end{Remark}

\subsection{Non expansiveness of Gabor frames}\label{s:8.2}

Consider now a Gabor frame $\G(\gamma,\alpha\Z^m\times \beta\Z^m)$,
 where $0<\alpha,\beta<1$ and $\gamma(x)=exp(-\norm{x}_2^2)$. The
 choice of $\al,\be$ will be irrelevant, but for the sake of example
 the reader may think to the case $\al=\be=\frac{1}{2}$. Let
 $\tilde{\gamma}$ denote its canonical dual frame generator. Let $E$ denote
the upper frame bound of $\G(\tilde{\gamma},\al\Z^m\times\be\Z^m)$. For two
 functions $f,h\in L^2(\R^{m})$, we denote by
\[ V_f h:\R^{2m}\rightarrow\C~~,~~ V_f h(\lambda)=\ip{h}{f_\lambda} \]
the windowed Fourier transform of $h$ with respect to $f$. The modulation
spaces $M^p$, $1\leq p\leq 2$, are defined by (see \cite{gr01}):
\[ M^p = \{ f\in L^2(\R^m)~|~V_\gamma f \in L^p(\R^{2m})\} ~~,~~\norm{f}_{M^p}
:= \norm{V_\gamma f}_{L^p} \]
In particular $\gamma,\tilde{\gamma}$ are both in $M^1$. 
Note $M^2=L^2$ as sets, and the norms are equivalent. 
The Wiener amalgam space $W(C,l^p)$ is defined by:
\[ W(C,l^p) = \{ f~;~f:\R^b\rightarrow\C~,~f~{\rm continuous}~,~
\norm{f}_{W(C,l^p)}^p := \sum_{k\in\Z^b}\sup_{x\in Q_1(k)}|f(x)|^p <\infty \} 
\]
The following result is proved in \cite{bacahela03-1}, Proposition A.3:
For all $f\in L^2(\R^m)$, $V_\gamma f\in W(C,l^2)$ and
\begin{equation}
\norm{V_\gamma f}_{W(C,l^2)}\leq C\norm{\gamma}_{M^1}\norm{f}_2
\end{equation}
where the constant $C$ can be chosen as $C=3^{m/2}$.
We can now prove the following.

\begin{Theorem}\label{t:gnonexp} Assume $\G(g,\Lambda)$ is a Gabor frame
in $L^2(\R^m)$. Then $\G(g,\Lambda)$ is non-expansive with respect to the
  quasi-metric index set $(\Lambda,\norm{\cdot}_\infty)$.
\end{Theorem}

\pf We will show the Gram operator of $\G$ is non-expansive, and then
the conclusion follows from Corollary \ref{c.c.1}.

We start with the following decomposition
\[ \ip{g_{\lambda_1}}{g_{\lambda_2}} = \sum_{k,j\in\Z^m}
\ip{g_{\lambda_1}}{\gamma_{\al k,\be j}}\ip{\tilde{\gamma}_{\al k,\be
j}}{ g_{\lambda_2}} = (AB)_{\lambda_1,\lambda_2} \] where
$A:l^2(\al\Z^{m}\times\be\Z^m)\rightarrow l^2(\Lambda)$, 
$B:l^2(\Lambda)\rightarrow l^2(\al\Z^{m}\times\be\Z^m)$, 
are defined through
$A_{\lambda,(\al k,\be j)}=\ip{g_\lambda}{\gamma_{\al k,\be j}}$,
$B_{(\al k,\be j),\lambda}=\ip{\tilde{\gamma}_{\al k,\be j}}{g_{\lambda}}$.
$A$ and $B$ are bounded operators since they are compositions of
analysis and synthesis operators associated to frames $\G(g,\Lambda)$,
$\G(\gamma,\al \Z^m\times \be\Z^m)$ and
$\G(\tilde{\gamma},\al\Z^m\times\be\Z^m)$. Note
\[ |A_{\lambda,(\al k,\be j)}| = |V_\gamma g((\al k,\be j)-\lambda)| \]
\[ |B_{(\al k,\be j),\lambda}| = |V_{\tilde{\gamma}} 
g(\lambda-(\al k,\be j))| \]

Consider the map $a:\Lambda\rightarrow\al\Z^{m}\times\be\Z^m$,
$a(\lambda)=(\zeta_k\lfloor\frac{\lambda_k}{\zeta_k}\rfloor)_{ 1\leq
k\leq 2m}$, where $\lambda=(\lambda_k)_{1\leq k\leq 2m}$,
$\zeta_k=\al$ for $1\leq k\leq m$, $\zeta_k=\be$ for $m+1\leq
k\leq 2m$, and $\lfloor x\rfloor$ is the largest integer smaller than
or equal to $x$. Thus $\norm{a(\lambda)- \lambda}_\infty < 1$.

Recall that $V_\gamma g$ and $V_{\tilde{\gamma}}g$ are both in
$W(C,l^2)$.  Combining this fact to the fact that every box of size
length $L_0$ has at most $U_0$ points (see previous subsection), we
obtain that, for every $\rho>0$ there are $N_A(\rho),N_B(\rho)>0$ so
that
\begin{equation}\label{eq:A}
 \forall r\in\al\Z^m\times\be\Z^m~~,~~
\sum_{\lambda\in\Lambda\setminus Q_{N_A}(r)} 
|V_\gamma g(r-\lambda)|^2 < \rho
\end{equation}
\begin{equation}\label{eq:B}
 \forall\lambda\in\Lambda~~,~~
\sum_{k,j\in\Z^m, \norm{(\al k,\be j)-a(\lambda)}_\infty >N_B(\rho)
}|V_{\tilde{\gamma}}g(\lambda - (\al k,\be j))|^2 < \rho
\end{equation}

Fix $\eps>0$. We will find
$N=N_\eps>0$ so that for all $\lambda\in\Lambda$,
\begin{equation}
\label{eq:desire}
 \sum_{\nu\in\Lambda\setminus B_N(\lambda)}|\ip{g_{\nu}}{g_\lambda}|^2 < 
\eps 
\end{equation}
Since the Gram operator is symmetric, this will conclude the proof.

The remainder of the proof mirrors the argument used in Theorem
\ref{t.c.1} that shows that non-expansiveness is preserved under
multiplication.

Let $\eps_B = \frac{\eps}{4\norm{A}^2}$ and $N_B=N_B(\eps_B)$ as in
(\ref{eq:B}),
$\eps_A=\frac{\eps}{4E\norm{g}^2}(\frac{\al\be}{(2N_B+1)^2})^m$ and
$N_A=N_A(\eps_A)$ the associated integer that satisfies (\ref{eq:A}).
Set $N=N_A+N_B+1$. We prove this choice satisfies (\ref{eq:desire}).

Let $(\delta_\lambda)_\lambda$ denote the sequence whose entries are
zero except for the $\lambda^{th}$ entry which is one. Thus
$\{\delta_\lambda~;~\lambda\in\Lambda\}$ is the canonical orthonormal
basis of $l^2(\Lambda)$. Note for all $\nu,\lambda\in\Lambda$,
$(AB)_{\lambda,\nu}=\ip{AB\delta_\nu}{\delta_\lambda}$. 

Fix a $\eta\in\Lambda$. Let $v,w\in l^2(\al\Z^m\times\be\Z^m)$ denote
the vectors of $B\delta_\eta=v+w$, where all entries of $v=(v_{\al
k,\be j})$ vanish for $\norm{(\al k,\be j)-a(\eta)}_\infty<N_B$, and
all entries of $w=(w_{\al k,\be j})$ vanish for $\norm{(\al k,\be
j)-a(\eta)}_\infty\geq N_B$. By (\ref{eq:B}) we obtain
$\norm{v}^2_{l^2}<\eps_B$, and hence
$\norm{Av}_{l^2}^2\leq\frac{\eps}{4}$. Now we have:
\[
 T:=\sum_{\lambda\in\Lambda\setminus B_N(\eta)}|(AB)_{\lambda,\eta}|^2  =
 \sum_{\lambda\in\Lambda\setminus B_N(\eta)}|\ip{Av}{\delta_\lambda}+
\sum_{\tiny \begin{array}{cc} \mbox{$r\in\al\Z^m\times\be\Z^m$} \\
\mbox{$\norm{r-a(\eta)}_\infty<N_B$} \end{array}}
A_{\lambda,r}B_{r,\eta}|^2 \]
\begin{eqnarray}
T & \leq & 2 \sum_{\lambda\in\Lambda}
|\ip{Av}{\delta_\lambda}|^2 + 2\sum_{\lambda\in\Lambda\setminus
B_N(\eta)}|\sum_{\tiny \begin{array}{cc} \mbox{$r\in\al\Z^m\times\be\Z^m$} \\
\mbox{$\norm{r-a(\eta)}_\infty<N_B$} \end{array}}
A_{\lambda,r}B_{r,\eta}|^2 \\
& \leq & \frac{\eps}{2} + 2 \sum_{\lambda\in\Lambda\setminus B_N(\eta)}
\left(\sum_{\tiny \begin{array}{cc} \mbox{$r\in\al\Z^m\times\be\Z^m$} \\
\mbox{$\norm{r-a(\eta)}_\infty<N_B$} \end{array}} 1\right)\left(
\sum_{\tiny \begin{array}{cc} \mbox{$r\in\al\Z^m\times\be\Z^m$} \\
\mbox{$\norm{r-a(\eta)}_\infty<N_B$} \end{array}}|A_{\lambda,r}B_{r,\eta}|^2
\right) \\
& \leq & \frac{\eps}{2} + 2(\frac{2N_B^2}{\al\be})^m 
\sum_{\tiny \begin{array}{cc} \mbox{$r\in\al\Z^m\times\be\Z^m$} \\
\mbox{$\norm{r-a(\eta)}_\infty<N_B$} \end{array}}|B_{r,\eta}|^2
\sum_{\lambda\in\Lambda\setminus B_N(\eta)}|A_{\lambda,r}|^2 \\
& \leq & \frac{\eps}{2} + 2\left(\frac{(2N_B+1)^2}{\al\be}\right)^m 
E\norm{g}^2 \eps_A = \frac{\eps}{2}+\frac{\eps}{2} = \eps
\end{eqnarray}
where the last inequality follows from $\Lambda\setminus B_N(\eta)\subset 
\Lambda\setminus Q_{N_A}(r)$, for all $r\in\al\Z^m\times\be\Z^m$ with
 $\norm{r-a(\eta)}_\infty<N_B$, and (\ref{eq:A}). This proves (\ref{eq:desire})
and thus the statement. \qed

\begin{Remark}
In terminology of \cite{bacahela06}, (\ref{eq:A}) means $(\G(g,\Lambda),a,
\G(\gamma,\al\Z^m\times\be\Z^m))$ has $l^2$-column decay, whereas 
(\ref{eq:B}) means that $(\G(g,\Lambda),a,
\G(\tilde{\gamma},\al\Z^m\times\be\Z^m))$ has $l^2$-row decay.

Using the terminology from \cite{bacahela06}, Theorem \ref{t:gnonexp}
states that $(\G(g,\Lambda),a)$ is $l^2$-self-localized, and
$l^2$-localized with respect to its canonical dual frame.
\end{Remark}

\subsection{Measure Functions of Gabor Superframes}\label{s:8.3}

Consider now two Gabor frames $\gc(g,\Lambda)$ and $\hc(h,\Sigma)$
in $L^2(\R^m)$. Assume there is a bijection $a:\Lambda\rightarrow\Sigma$
so that $(\gc(g,\Lambda),\hc(h,\Sigma))$ is a superframe, that is
\begin{equation}\label{eq:9.3.1}
\fc = \{g_\lambda \oplus h_{a(\lambda)}~;~\lambda\in\Lambda\}
\end{equation}
is frame for $L^2(\R^m)\oplus L^2(\R^m)$. Note $\fc\in\fc[\Lambda]$.

\begin{Proposition}
\label{t9.3.1}
Assume $(\gc(g,\Lambda),\hc(h,\Sigma))$ is a Gabor superframe with respect
to the correspondence $a:\Lambda\rightarrow\Sigma$. Assume there are
invertible matrices $M^1,M^2\in\R^{2m\times 2m}$so that the map $a$
satisfies
\begin{equation}
\label{eq:9.3.22}
\lim_n\frac{|a^{-1}(\Sigma\cap(M^2Q_n(O)))\cap(M^1Q_n(O))|}{|\Lambda
\cap(M^1Q_n(O))|} = \lim_n \frac{|\Sigma\cap(M^2 Q_n(O))|}{|\Lambda
\cap(M^1 Q_n(O))|}=1
\end{equation}
and there exists a function $r:[0,\infty) \rightarrow [0, \infty)$ such that 
for all $\sigma_1,\sigma_2\in\Sigma$,
\begin{equation}
\label{eq:9.3.23}
\norm{a^{-1}(\sigma_1)-a^{-1}(\sigma_2)}\leq 
r(\norm{\sigma_1-\sigma_2})
\end{equation}
Then the direct sum frame $\fc$ defined in (\ref{eq:9.3.1}) has 
the free ultrafilter frame measure:
\begin{equation}
\label{eq:9.3.2}
\mu(\fc)(p) = \frac{1}{D(\Lambda;p,M^1)}+\frac{1}{D(\Sigma;p,M^2)}
~~~,~~~\forall p\in\N^*
\end{equation}
In particular, the following is a necessary condition:
\begin{equation}
\label{eq:9.3.24}
limsup_n (\frac{det(M^1)}{|\Lambda\cap (M^1Q_n(O))|}+\frac{det(M^2)}{
|\Sigma\cap(M^2Q_n(O))|})n^{2m} \leq 1
\end{equation}
\end{Proposition}

\pf Note (\ref{eq:9.3.22}) and (\ref{eq:9.3.23}) imply that $a$ satisfies 
(\ref{eq:8.2.1}) and (\ref{eq:8.3.1}). Now (\ref{eq:9.3.2}) follows from 
Theorems \ref{t8.1}, \ref{t:gabormeasure}, and \ref{t:gnonexp}. 
Equation (\ref{eq:9.3.24})
 is obtained from (\ref{eq:9.3.2}), and (\ref{eq:dens}),
 and the fact that for any frame $\fc$, $\mu(\fc)(p)\leq 1$ for all $p$.
 \qed

\begin{Remark} 
Let $L_\Sigma>0$ be such that any box of side length
 $L_\Sigma$ in $\R^{2m}$ contains at least one point of $\Sigma$. 
Then condition (\ref{eq:9.3.23}) can be replaced equivalently by the following
boundedness condition:
\begin{equation}
\label{eq:9.3.25}
\exists R_0>0,\forall \sigma_1,\sigma_2\in\Sigma~~
\norm{\sigma_1-\sigma_2}\leq \sqrt{2m}L_\Sigma~\Rightarrow~
\norm{a^{-1}(\sigma_1)-a^{-1}(\sigma_2)}\leq R_0
\end{equation}
Indeed, if (\ref{eq:9.3.25}) holds true then for any $N>0$ there is a
chain of $\frac{N}{\sqrt{2m}L_\Sigma}$ points in $\Sigma$ so that the distance
between any two adjacent points is at most $\sqrt{2m}L_\Sigma$. Using the
 triangle
inequality it follows that (\ref{eq:9.3.23}) is satisfied with
$r(u)=(1+\frac{u}{\sqrt{2m}L_\Sigma})R_0$.
\end{Remark}

Using induction  one can immediately prove:
\begin{Theorem}
\label{t9.3.2}
Assume $\gc(g^k,\Lambda_k)$, $1\leq k\leq d$, are Gabor frames in
$L^2(\R^m)$ so that for maps $a_k:\Lambda_1\rightarrow\Lambda_k$,
$2\leq k\leq d$, the set $\fc=\{g^1_\lambda\oplus
g^2_{a_2(\lambda)}\oplus\cdots\oplus
g^d_{a_d(\lambda)}~;~\lambda\in\Lambda_1\}$ is frame for
$L^2(\R^m)\oplus\cdots\oplus L^2(\R^m)$. Assume further that there are
invertible matrices $M^k$, $1\leq k\leq d$ such that all maps $a_k$
satisfy
\begin{equation}
\label{eq:9.3.26}
\lim_n\frac{|a_k^{-1}(\Lambda_k\cap(M^kQ_n(O)))\cap(M^1Q_n(O))|}{|\Lambda_1
\cap(M^1Q_n(O))|} = \lim_n \frac{|\Lambda_k\cap(M^k Q_n(O))|}{|\Lambda_1
\cap(M^1 Q_n(O))|}=1
\end{equation}
and there exists a map $r:[0,\infty) \rightarrow [0,\infty)$ such that 
for all $\sigma_1,\sigma_2\in\Lambda_k$,
\begin{equation}
\label{eq:9.3.27}
\norm{a_k^{-1}(\sigma_1)-a_k^{-1}(\sigma_2)} \leq 
r(\norm{\sigma_1-\sigma_2})
\end{equation}
Then the free ultrafilter frame measure function of $\fc$ is given by
\begin{equation}
\label{eq:9.3.4}
\mu(\fc)(p)=\frac{1}{D(\Lambda_1;p,M^1)}+\cdots+\frac{1}{D(\Lambda_d;p,M^d)}~~,~~
p\in\N^*
\end{equation}
In particular it follows that necessarily
\begin{equation}
\label{eq:9.3.5}
\frac{1}{D(\Lambda_1;p,M^1)}+\cdots+\frac{1}{D(\Lambda_d;p,M^d)} \leq 1~~,~~
\forall p\in\N^*
\end{equation}
\end{Theorem}

In the special case of regular Gabor frames, $\Lambda_k=\{A_k n~;~
n\in\Z^{2m}\}$, $1\leq k\leq d$, we obtain that if
$(\gc(g_1;\Lambda_1),\ldots,\gc(g_d;\Lambda_d))$ form a superframe
with respect to the maps $a_k:\Lambda_1\rightarrow\Lambda_k$, $a_k(A_1n)=A_kn$,
 $2\leq k\leq d$, then conditions (\ref{eq:9.3.26}) and (\ref{eq:9.3.27})
are satisfied with $M^k=A_k$,  and we obtain
immediately the following result which recovers and extends the result
of \cite{ba99-1},

\begin{Corollary}\label{c9.4}
Assume $g^1,\ldots,g^d\in L^2(\R^m)$ and $A_1,\ldots,A_d\in\R^{2m\times 2m}$
are so that $\fc=\gc(g^1,A_1\Z^{2m})\oplus\cdots\oplus\gc(g^d,A_d\Z^{2m})$
is frame for $L^2(\R^{m})\oplus\cdots\oplus L^2(\R^m)$, then for any frame
measure function $m:\fc[\Z^{2m}]\rightarrow\cstarw$,
\begin{equation}
\label{eq:9.4.1}
m(\fc) = (det(A_1)+\cdots+det(A_d))1_{W}
\end{equation}
Consequently, as a necessary condition to have a superframe,
\begin{equation}
\label{eq:9.4.2}
det(A_1)+\cdots+det(A_d)\leq 1
\end{equation}
\end{Corollary}

\section{Redundancy} {\label{s:red}}

The word {\it redundancy} is often used to describe, qualitatively, the
overcompleteness of frames.  However, for frames with an infinite
number of elements, there is no quantitative definition of redundancy.
Here, we propose that the reciprocal of a frame measure function
should be the quantitative definition of redundancy.

\begin{Definition}
Given a measure function $m:\fci \rightarrow \cstarm$, we define the
{\em redundancy function} $R:\fci \rightarrow \mbox{\{functions from
$M$ to $\R \cup \infty$\}}$, $R(\fc)(x)= (m(\fc)(x))^{-1}$.  In the
case when the measure function is the ultrafilter measure function, we
term the redundancy function the {\em ultrafilter redundancy
function.}
\end{Definition}

The rest of this section discusses the justification for this
definition.  We begin by listing a series of properties of the frame
redundancy function, all of which mesh well with the qualitative
notion of redundancy:
\begin{itemize}
\item We immediately have the desirable properties that for a frame,
the redundancy function is greater than or equal to one with the
redundancy function equal to one for any Riesz basis.

\item By Theorem \ref{t:gabormeasure}, for any Gabor frame
$\G(g,\Lambda)$ and indexing $(\Lambda,(\Lambda_n)_n)$ as in Section 
\ref{s:gab}, the
ultrafilter redundancy function corresponds to the density of the time
frequency shifts as follows:
\begin{equation}
\label{eq:9.1.1}
 R(\G(g,\Lambda)(p)= D(\Lambda; p,M), 
\mbox{ for all free ultrafilters } p.
\end{equation}

\item This connection between redundancy and measure function extends
to localized frames.  Using the notation and results from \cite{bacahela06}
we have an explicit description of the ultrafilter redundancy
function. Assume $\fc\in\fc[I]$ is a frame for $H$ and
$a:I\rightarrow\Z^d$ is a map so that $(\fc,a,\ec)$ has both
$l^2$-column and $l^2$-row decay (see \cite{bacahela06} for definition), 
where $\ec=\{e_k~;~k\in\Z^d\}$ is another frame for $H$.  
Set $I_n=a^{-1}(Q_n(0))$, where $Q_n(0)$ is
the box of side length $n$ centered at 0 in $\Z^d$, and consider the
ultrafilter redundancy functions associated to $(I,(I_n)_n)$,
respectively $(Z^d,(Q_n(0))_n)$.  Then Theorem 5 in \cite{bacahela06}
implies:
\begin{equation}
\label{eq:9.1.2}
R(\fc)(p) = D(a;p)R(\ec)(p)
\end{equation}
In particular, if $\ec$ is a Riesz basis for $H$, then $R(\ec)=1$ and
the previous equation turns simply into:
\begin{equation}
\label{eq:9.1.3}
R(\fc)(p) = D(a;p)
\end{equation}

\item In these cases (Gabor and localized frames), the redundancy
function is additive on unions of frames. Suppose $\fc^1\in\fc[I]$
and $\fc^2\in\fc[J]$ are two frames for same Hilbert space $H$, and that
there are maps $a^1:I\rightarrow\Z^d$ and $a^2:J\rightarrow\Z^d$ so
that $(\fc^1,a^1,\ec)$ and $(\fc^2,a^2,\ec)$ have both $l^2$-column
and $l^2$-row decay, where $\ec$ is a Riesz basis for $H$. Set
$I_n=(a^1)^{-1}(Q_n(0))$, and $J_n=(a^2)^{-1}(Q_n(0))$.  Consider the
ultrafilter redundancy functions associated to $(I,(I_n)_n)$,
$(J,(J_n)_n)$, $(I\cup J,(I_n\dot{\cup} J_n)_n)$ for frames $\fc^1$,
$\fc^2$, and $\fc^1\dot{\cup}\fc^2$, respectively. Here $\dot{\cup}$
denotes union with multiplicity. 
First it is immediate to check that
$(\fc^1\dot{\cup}\fc^2,a,\ec)$ has $l^2$-column and
$l^2$-row decay, where $a:I\dot{\cup}J\rightarrow\Z^d$,
$a(i)=a^1(i)$ for $i\in I$, and $a(j)=a^2(j)$ for $j\in J$. 
Next note that $I_n\dot{\cup} J_n=a^{-1}(Q_n(0))$. Then, applying
(\ref{eq:9.1.3}) to $\fc^1\dot{\cup}\fc^2$ we obtain:
\begin{equation}
\label{eq:9.1.4}
R(\fc^1\dot{\cup}\fc^2)(p)= D(a;p)=D(a^1;p)+D(a^2;p) = R(\fc^1)(p)+
R(\fc^2)(p)
\end{equation}
which proves additivity of the redundancy function. Equation (\ref{eq:9.1.4})
can be immediately extended to any finite number of frames.
\end{itemize}

In addition to the above properties, the redundancy function can be
seen as an analogue of redundancy in the finite dimensional case.  In
finite dimensions, the idea of redundancy is quantified.  Here we have
a frame $\fc=\{f_j\}_{j\in J}$, consisting of $M=|J|$ vectors.  If we
let $N$ be the dimension of the space spanned by the elements of
$\fc$, then the ratio $r=\frac{M}{N}$ is a natural quantity that is
often referred to as the redundancy of the frame $\fc$.  Another way
to arrive at the quantity $r$ is as follows.  Associated to $\fc$ is
the finite dimensional Gram operator $G:l^2(J) \rightarrow l^2(J)$
defined entry-wise by $G_{i,j}=\ip{f_i}{f_j}$. The ratio of the
dimension of the space $l^2(J)$ (which is $|J|$) to the dimension of
the range of $G$ is also $r=\frac{M}{N}$.  In other words the
reciprocal of the redundancy, $\frac{1}{r}$, is the normalized trace
of the associated Gram projection of the frame.

So what is the meaning of the quantity $r$?  In this setting we have
that a frame $\fc$ is a basis if and only if $r=1$.  If $\fc$ is the
union of two bases on the same space then $r=2$, however this is not
the only type of frame that has $r=2$; a basis of size $n$ along with
$n$ additional copies of the first basis element also has $r=2$.  Thus
the value of $r$ does not reveal the whole story, but it does provide
a one paramater classification of frames.  One can then examine the
set of frames with a given $r$ and try and understand the variation in
their characteristics (see \cite{befi03,cako03}).  One can also design frames
with a particular value of $r$ that maximizes certain channel capacity
or energy considerations \cite{hest03,tdhs05}.

If one tries to use the finite dimensional case as a road map for
defining redundancy in infinite dimensions, one immediately encounters
difficulty.  In this case, we are considering a frame $\fc = \{f_i
\}_{i\in I}$ indexed by an infinite set $I$.  Thus the corresponding
quantity $M=|I|$ is infinite.  Generically, the dimension of the space
spanned by the $f_i$ which was denoted by $N$ in the finite case is
also infinite and therefore the ratio $r=\frac{M}{N}$ is meaningless.
Similarly, attempting to compare the dimension of $l^2(I)$ to the
dimension of the range of the Gram operator of $\fc$, yields a
comparison of two infinite quantities.

By itself, comparing the dimensions of infinite dimensional spaces is
not completely hopeless.  Those familiar with the study of von Neumann
algebras will recall that the dimension function, introduced by von
Neumann, provides a way of comparing certain infinite dimensional
subspaces of a fixed infinite dimensional space.  In this case, only
subspaces that are ranges of projections in the algebra are
considered; the dimension function of the subspace is then defined to
be the normalized trace (which exists on a Von Neumann algebra) of the
projection.  This connection has yielded many nice results about Gabor
frames on regular lattices \cite{ri81-1,dalala95,ja95,feko98,feka04} (just
to name a few); in these cases the regular lattice structure was
enough to ensure that the Gramian had the necessary structure to allow
the tools of von Neumann algebras to be useful.  In general, however,
this added structure is not available and we are further discouraged
by the known fact that there does not exist a dimension function that
is finite and non-zero on all non-zero subspaces of a fixed infinite
dimensional space.

As mentioned earlier, in finite dimensions the reciprocal of the
redundancy can be defined as the trace of the associated Gram
projection to the given frame.  The ultrafilter redundancy function
can be seen as the infinite dimensional analogue of this.  To begin
with, the ultrafilter frame measure function is determined by certain
averages of $\ip{\tilde{f}_i}{f_i}$, that is, certain averages of the
diagonal elements of the corresponding Gram projection-- a natural
generalization of the normalized trace in finite dimensions which is
the average of the diagonal elements of the Gram projection.  The key
structural feature of a trace is that the trace of $AB$ and $BA$ are
equal for operators $A$ and $B$.  This feature is present for measure
functions on the set of non-expansive operators (Lemma \ref{l.c.4}).

For these reasons, we feel our definition is the proper quantification
of redundancy in the infinite setting.  There remain unanswered
questions about the redundancy function, an important one being if a
frame has redundancy $c$, does there exist a subset of the frame that
is a frame for the same space with redundancy $1$ (or $1 + \epsilon$
for any $\eps>0$).

\appendix
\section{Supersets} \label{a:supersets}

We recall the notion of superframe (see \cite{ba98-1,ba99-1,ba00})
(or disjoint frames, as used by D.Larson, see \cite{hala00}). 
Let $\fc_1,\ldots,\fc_L\in\fci$, a finite number of frames indexed by $I$.

\begin{Definition}
We call $(\fc_1,\ldots,\fc_L)$ a {\em superframe} if
\begin{equation}
\label{eq.m1.1}
\fc = \fc_1\oplus\cdots\oplus\fc_L:=\{ f^1_i\oplus\cdots\oplus f^L_i~;~i\in
I\}
\end{equation}
is a frame in $H_1\oplus\cdots\oplus H_L$, the direct sum of Hilbert spaces
spanned by $\fc_1,\ldots,\fc_L$, respectively.
\end{Definition}

An equivalent characterization of superframes is given by the following
\begin{Theorem}[\cite{ba00}]\label{p:A2}
The collection $(\fc_1,\ldots,\fc_L)$ is a superframe if and only if the
following two conditions hold true:
\begin{enumerate}
\item Each $\fc_l$ is frame, $1\leq l\leq L$;
\item $E_k\cap(\sum_{l\neq l}E_l)=\{0\}$, for $1\leq k\leq L$, and
$\sum_{k=1}^LE_l$ is closed (where $E_l$ is the range in $\ltwoi$ 
of the analysis operator associated to $\fc_l$).
\end{enumerate}
\end{Theorem}

In particular, the second condition above holds true when the 
ranges of $E_l$ are mutually orthogonal. This special case is called {\em
orthogonal in the sense of supersets} (or strongly disjoint, see
\cite{hala00}). More specifically we define the following:
\begin{Definition}
Two frames $\fc_1=\{f^i_i;i\in I\}$ and $\fc_2=\{f^2_i;i\in
I\}$ indexed by $I$ are said to be {\em orthogonal in the sense of
supersets} if $E_1$, the range of analysis operator associated to
$\fc_1$, is orthogonal in $\ltwoi$ to $E_2$, the range of coefficients
associated to $\fc_2$.  Equivalently,
\begin{equation}
\label{eq.m1.2}
\sum_{i\in I}\ip{g}{f^1_i}\ip{f^2_i}{h} =0~~,~~\forall g\in H_1~,~\forall h\in
H_2
\end{equation}

\end{Definition}
\begin{Remark} \label{r.m1.1}
Clearly if two frames $\fc_1,\fc_2$ are orthogonal in the sense of
supersets, then $E_1\cap E_s=\{0\}$ and $E_1+E_2$ is closed, hence
$(\fc_1,\fc_2)$ is a superframe. Note that in this case the range of
the analysis operator associated to $\fc_1\oplus\fc_2$ is exactly
$E_1\oplus E_2$, and the associated Gram projection $P$, is given by
$P=P_1+P_2$, the sum of the associated Gram projections of $\fc_1$ and
$\fc_2$. In particular, the canonical dual of $\fc_1\oplus\fc_2$ is the
direct sum of the canonical duals of $\fc_1$ and $\fc_2$.
\end{Remark}
\begin{Remark}
\label{r.m1.2}
For any frame $\fc\in\fci$, one can always construct $\fc'\in\fci$
that is orthogonal to $\fc$ in the sense of supersets. Let $P$ be the
associated Gram projection to $\fc$. Then $Q=1-P$ is also an
orthogonal projection in $\ltwoi$ (1 being the identity operator).
Set $\fc'=\{Q\delta_i~;~i\in I\}$. One can easily check that $\fc'$ is
a (Parseval) frame and that its associated Gram projection is $Q$;
therefore $\fc$ and $\fc'$ are orthogonal in the sense of supersets.
\end{Remark}

\section{Ultrafilters} \label{a:ultr}

Consider the difference between the limit of a sequence and the liminf
of a sequence.  The liminf has the advantage that it is defined on all
bounded sequences as opposed to the limit which is only defined on the
relatively small set of sequences that have limits.  However, unlike
the limit, the liminf is not linear on its domain.

The existence of ultrafilters leads to linear functionals (Definition
\ref{d:uf} that achieve ``the best of both worlds'' in the sense that
they are defined and linear on all bounded sequences (Proposition
\ref{p.m2.2}).

\begin{Definition}
A collection $\accc$ of subsets of $M$ is called
a {\em filter} if it satisfies the following properties:
\begin{enumerate}
\item The empty set is not in $\accc$: $\emptyset\not\in\accc$;
\item If $A_1,A_2\in\accc$, then $A_1\cap A_2\in\accc$;
\item If $A \subset B \subset M$ with $A \in \accc$ then $B \in \accc$.
\end{enumerate}

A filter $\accc$ is an {\em ultrafilter} if it is 'maximal' in the
following sense:
\begin{enumerate}
\setcounter{enumi}{3}
\item 
 For all $A \subset M$ either $A \in\accc$ or $(M\setminus A)\in\accc$ 
(but not both  because of 1. and 2. above).
\end{enumerate}
An ultrafilter that does not contain a finite set is called a {\it
free} ultrafilter; the set of free ultrafilters will be denoted by
$M^*$.
\end{Definition}

The existence of free ultrafilters is unintuitive and requires the
axiom of choice.  For our purposes we shall be concerned with the case
$M=\N$, and $\N^*$ denotes the set of free ultrafilters.

The existence of ultrafilters allows us to define a family of limits
on bounded sequences indexed by $M$:
\begin{Definition} \label{d:uf}
Let $\xv=\{x_m \}_{m\in M}$ be a  bounded sequence of complex numbers.   
Given an ultrafilter $p$ on $M$, we say {\em $x$ converges to $c\in\C$
with respect to the ultrafilter $p$} and write $c=\plim\ \xv$,
if for any $\eps>0$ there is a set $A\in\accc$ such that $|x_m-c|<\eps$ 
for all $m\in A$. 
\end{Definition}

This notion of limit has the following consequences that can be found
in any text about ultrafilters (see \cite{hindeman98} for example):

\begin{Proposition}\label{p.m2.2}
Let $\xv=\{x_m \}_{m\in M}$, $\yv=\{y_m \}_{m\in M}$ be bounded
sequences of complex numbers and let $p$ be a free ultrafilter.
\begin{enumerate}
\item   $\plim\ \xv$ exists and is unique.
\item The function $\plim$ is linear, i.e. $\plim (a\xv + b\yv)=
a(\plim\ \xv)) + b (\plim\ \yv)$ for all scalars $a,b$.
\item For $M=\N$, the value of $\plim\ \xv$ is an accumulation point
of the set $x_1, x_2, \dots$ Consequently, if the sequence $x_1, x_2,
\dots$ has a limit, then $\plim\ \xv$ is equal to that limit.
\end{enumerate}
\end{Proposition}

\section{Acknowledgments}

The authors acknowledge very useful discussions with Pete Casazza, Ingrid
Daubechies, Hans Feichtinger, Sinan Gunturk, Christopher Heil, 
Gitta Kutyniok, and Henry Landau.

\bibliographystyle{plain}
\bibliography{nhgbib2}

\end{document}